\documentclass [11pt,a4paper]{article}
\usepackage[ansinew]{inputenc}
\usepackage{amssymb}
\usepackage{amsmath}
\usepackage{graphicx}
\usepackage{amsthm}
\usepackage{cite}

\newtheorem{theorem}{\textbf{Theorem}}[section]
\newtheorem{lemma}{\textbf{Lemma}}[section]
\newtheorem{proposition}{\textbf{Proposition}}[section]
\newtheorem{corollary}{\textbf{Corollary}}[section]
\newtheorem{remark}{\textbf{Remark}}[section]
\newtheorem{definition}{\textbf{Definition}}[section]

\allowdisplaybreaks[4]

\def\be{\begin{equation}}
\def\ee{\end{equation}}
\def\bea{\begin{eqnarray}}
\def\eea{\end{eqnarray}}
\def\bt{\begin{theorem}}
\def\et{\end{theorem}}
\def\bl{\begin{lemma}}
\def\el{\end{lemma}}
\def\br{\begin{remark}}
\def\er{\end{remark}}
\def\bp{\begin{proposition}}
\def\ep{\end{proposition}}
\def\bc{\begin{corollary}}
\def\ec{\end{corollary}}
\def\bd{\begin{definition}}
\def\ed{\end{definition}}
\def\non{\nonumber }

\begin{document}

\title{Long-time dynamics of the nonhomogeneous incompressible flow of nematic liquid crystals}

\author{{\sc Xianpeng Hu}\footnote{Courant Institute of Mathematical Sciences, New York University, New York, NY 10012, USA.
E-mail: \textit{xianpeng@cims.nyu.edu}} \ 
and {\sc Hao Wu}\footnote{Corresponding author. School of
Mathematical
Sciences and Shanghai Key Laboratory for Contemporary Applied Mathematics, Fudan University, Shanghai 200433, China. E-mail: \textit{haowufd@yahoo.com}}}

\date{\today}

\maketitle


\begin{abstract}
We study the long-time behavior of global strong solutions to a hydrodynamic system for nonhomogeneous incompressible nematic liquid crystal flows driven by two types of external forces in a smooth bounded domain in $\mathbb{R}^2$. For arbitrary large regular initial data with the initial density being away from vacuum, we prove the decay of the velocity field for both cases. Furthermore, for the case with asymptotically autonomous external force, we can prove the convergence of the density function and the director vector as time goes to
infinity. Estimates on convergence rate are also provided.


\end{abstract}

\section{Introduction}
Liquid crystals are substances that exhibit a phase of matter that has properties between
those of a conventional liquid, and those of a solid crystal \cite{dG}. The hydrodynamic theory of liquid crystals
due to Ericken and Leslie was developed around 1960's
\cite{E61,Le68,Le79}. Since then, the mathematical theory is still progressing
and the study of the full Ericksen--Leslie model presents relevant mathematical difficulties. We consider the following hydrodynamical model for the flow of
nematic liquid crystals (cf. \cite{lin1})
 \bea
 \rho_t+v \cdot \nabla \rho&=&0,\label{0}\\
 \rho(v_t+v\cdot\nabla v)-\nu \Delta v+\nabla P&=&-\lambda
 \nabla\cdot(\nabla d\odot\nabla d)+\rho\mathbf{g},\label{1}\\
 \nabla \cdot v &=& 0,\label{2}\\
 d_t+v\cdot\nabla d&=&\gamma(\Delta d-f(d)),\label{3}
 \eea
in $\Omega \times\mathbb{R}^+$, where $\Omega \subset \mathbb{R}^n$ $(n=2,3)$ is assumed to be a bounded domain with smooth boundary $\Gamma$. System \eqref{0}--\eqref{3} is subject
to the Dirichlet boundary conditions:
 \be
 v(x,t)=0,\quad d(x,t)=d_0(x),\qquad \text{for}\ (x, t)\in \Gamma\times
 \mathbb{R}^+.
 \label{4}
 \ee
and the initial conditions
 \be
 \rho|_{t=0}=\rho_0(x), \ \ \  v|_{t=0}=v_0(x) \ \ \text{with}\ \nabla\cdot v_0=0,\quad
 d|_{t=0}=d_0(x),\qquad \text{for}\ x\in \Omega.\label{5}
 \ee
In the above system, $\rho$ is the density of the material,
$v$ is the velocity field of the flow and $d$ represents the
averaged macroscopic/continuum molecular orientation in
$\mathbb{R}^n$. $P(x,t)$ is a scalar function representing
the pressure (including both the hydrostatic and the induced elastic
part from the orientation field). $\mathbf{g}$ stands for the external body force. The positive constants $\nu,
\lambda$ and $\gamma$ stand for viscosity, the competition between
kinetic energy and potential energy, and macroscopic elastic
relaxation time (Debroah number) for the molecular orientation
field. We assume that $f(d)=\nabla F(d)$ for some smooth bounded
function $F:\mathbb{R}^n\rightarrow \mathbb{R} $. $\nabla d\odot
\nabla d$ denotes the $n\times n$ matrix whose $(i,j)$-th entry is
given by $\nabla_i d\cdot \nabla_j d$, for $1\leq i,j\leq n$.

System \eqref{0}--\eqref{3} is a nonhomogeneous version of the following simplified system introduced in \cite{lin1, LL95} that
models the incompressible flow of nematic liquid crystals with varying director lengths:
 \bea
 v_t+v\cdot\nabla v-\nu \Delta v+\nabla P&=&-\lambda
 \nabla\cdot(\nabla d\odot\nabla d),\label{a1}\\
 \nabla \cdot v &=& 0,\label{a2}\\
 d_t+v\cdot\nabla d&=&\gamma(\Delta d-f(d)).\label{a3}
 \eea
System \eqref{a1}--\eqref{a3} keeps the important mathematical
structure as well as some of the essential features of the
original Ericksen--Leslie system. The Ginzburg--Landau approximation
 \be
f(d)=\frac{1}{\eta^2}(|d|^2-1)d,\ \text{with its
antiderivative}\  F(d)=\frac{1}{4\eta^2}(|d|^2-1)^2,   \label{FF}
 \ee
 was introduced in order to relax the nonlinear constraint
$|d|=1$. The system \eqref{a1}--\eqref{a3} has been
studied in a series of work not only theoretically \cite{LL95,LL96, LS01,W10} but also numerically \cite{Liu2,Liu1} (see also \cite{HW10} for the case $f(d)=0$). In \cite{LL95}, the authors proved the existence of global weak solutions
to system \eqref{a1}--\eqref{a3} with Dirichlet boundary conditions by a semi-Galerkin method. Global existence and
uniqueness of classical solutions to the same system was proved for $n = 2$ or
$n = 3$ under large viscosity assumption. Long-time behavior of global solutions to system \eqref{a1}--\eqref{a3} was studied in \cite{LL95, PRS,W10}. In particular, convergence of
global classical solutions to single steady states as time goes to
infinity was obtained in \cite{W10,PRS}.  We refer to \cite{B,LS01, CGM09,GW12} for results on the homogeneous system \eqref{a1}--\eqref{a3} subject to other types of boundary conditions.

As far as the density-dependent system \eqref{0}--\eqref{3} is concerned, the authors in \cite{LZ09, Tan,XTan11} proved existence of global weak solutions of the problem \eqref{0}--\eqref{5} without assuming the positive lower bound for the initial density. The basic idea of their proof is to introduce a viscous term $\epsilon\Delta \rho$ in the transport equation \eqref{0}, and then pass to the limit as $\epsilon\to 0$. In the recent work \cite{DQS}, instead of introducing the viscosity term in \eqref{0}, the authors provided an alternative proof for the existence of global weak solutions to system \eqref{0}--\eqref{5} under a stronger assumption \eqref{bd} (i.e., the initial density is positive and bounded). Regularity properties of weak solutions to system \eqref{0}--\eqref{5} was proved by using Ladyzhenskaya type energy estimates for the approximate solutions constructed within a proper Galerkin scheme, provided that the initial data are regular and satisfy assumptions \eqref{bd}--\eqref{111}. For compressible version of the liquid crystal system \eqref{0}--\eqref{5}, existence and large-time behavior of a global weak solution were established in \cite{CML12, LQ12, WDYC12} while existence of local strong solutions was obtained in \cite{LLH11} (see also \cite{LXGLNM} for a blow-up criterion). Finally, we refer to the recent work \cite{CLL12, HuangWang12, LXWD12, LLW, LW10, Wang11, WD11, XZ12} and the references cited therein for mathematical results on the liquid crystal system under constraint $|d|=1$.

 We note that the external force $\rho\mathbf{g}$ is supposed to be vanishing in the above mentioned work. In this paper, we focus on the two-dimensional case $n=2$ and extend the results on long-time behavior of global classical solutions in \cite{LL95, W10} to the nonhomogeneous system \eqref{0}--\eqref{3} with non-vanishing external forces. Two types of external forces will be treated in the following text:
\begin{itemize}
\item[(F1)] $\mathbf{g}$ is a time-independent potential field, namely,
\be \mathbf{g}=\nabla \phi, \quad \text{for some scalar function}\ \phi(x)\in H^2(\Omega).\non\ee
\item[(F2)] $\mathbf{g}$ depends on time and satisfies the following integrability conditions:
 \bea
 && \mathbf{g}\in L^2(0,+\infty; \mathbf{H}^1(\Omega)), \quad \mathbf{g}_t\in L^2(0, +\infty; \mathbf{L}^2(\Omega)).\non
  \eea
\end{itemize}

The main results of this paper are as follows:
\bt
 \label{main2d}
 Suppose that $r\in (2,+\infty)$, the external force $\mathbf{g}$ satisfies either (F1) or (F2).
  For any initial data $\rho_0\in W^{1,r}(\Omega)$, $v_0\in \mathbf{H}^2(\Omega)\cap V$ and $d_0\in \mathbf{H}^3(\Omega)$ that satisfy
 \bea
 &&0<\underline{\rho}\leq \rho_0(x)\leq \bar{\rho},\quad \forall\,
 x\in \Omega,\label{bd}\\
 &&|d_0(x)|\leq 1,\quad \forall\, x\in \Omega,\label{111}
 \eea
 where $\underline{\rho}$ and $\bar{\rho}$ are positive constants, problem \eqref{0}--\eqref{5} admits a unique global strong solution $(\rho, v, d)$ such that for any $T>0$
 \bea
  && \rho \in C([0,T], W^{1,r}(\Omega)), \non\\
  && v\in C([0,T]; \mathbf{H}^2(\Omega)\cap V)\cap L^2(0,T; \mathbf{H}^3), \quad v_t\in L^2(0,T; V), \non\\
  && d\in C([0,T]; \mathbf{H}^3(\Omega))\cap L^2(0,T; \mathbf{H}^4), \quad d_t\in C([0,T]; \mathbf{H}^1_0(\Omega))\cap L^2(0,T; \mathbf{H}^2),\non\\
  && 0<\underline{\rho}\leq \rho(x,t )\leq \bar{\rho},\quad |d(x,t)|\leq 1,\quad \forall\,
 (x,t)\in \Omega\times[0,T].\non
 \eea
 \et

 \bt
 \label{conv1}
 Suppose that the assumptions of Theorem \ref{main2d} are satisfied with $\mathbf{g}$ fulfilling (F1).
  The global strong solution to problem \eqref{0}--\eqref{5} has the
following property
 \be \lim_{t\rightarrow +\infty}
 (\|v(t)\|_{\mathbf{H}^1}+\|v_t\|+\|d_t\|_{\mathbf{H}^1})=0.\label{conva}
 \ee
 For any unbounded sequence $\{t_i\}$, there is a subsequence $\{t'_i\}\nearrow+\infty$ such that
 \bea
 && \|\rho(t'_i)- \rho_\infty\|_{L^q}\to 0, \quad  \text{as}\ t'_i\to+\infty, \ q\in (1,+\infty), \label{conrho1}\\
 && \|d(t'_i)-d_\infty\|_{\mathbf{H}^3}\to 0,\quad \text{as}\ t'_i\to+\infty,\label{condr1}
 \eea
 where $\rho_\infty$ is a certain function that belongs to $L^q$ and $d_\infty$ is a solution to the following nonlinear elliptic boundary value
  problem:
  \be \left\{\begin{array}{l}  - \Delta d_\infty + f(d_\infty)=0,\quad x\in \Omega, \\
   d_\infty=d_0(x), \;\;\; x\in \Gamma.\\
   \end{array}
   \label{staa}
  \right.
 \ee
  \et

 \bt
 \label{conv}
 Suppose that the assumptions of Theorem \ref{main2d} are satisfied with $\mathbf{g}$ fulfilling (F2). If in addition, $\mathbf{g}$ satisfies
 \be
 \sup_{t\geq 0} (1+t)^{1+\xi} \int_t^{+\infty}\|\mathbf{g}(\tau)\|^2 d\tau<+\infty, \quad \text{for some}\ \xi>0,\label{ggg}
 \ee
 then the global strong solution to problem \eqref{0}--\eqref{5} has the
following property
 \be \lim_{t\rightarrow +\infty}
 (\|\rho(t)-\rho_\infty\|_{L^q}+\|v(t)\|_{\mathbf{H}^1}+\|v_t\|+\|d_t\|_{\mathbf{H}^1}+\|d(t)-d_\infty\|_{\mathbf{H}^3})=0,\label{cgce}
 \ee
 where $\rho_\infty$ is certain function in $L^q$ $(q\in (1,+\infty))$ and $d_\infty$ is a solution of \eqref{staa}.
 Moreover, there exists a positive constant $C$ depending on
 $v_0,d_0, \nu, \eta, \Omega, \bar{\rho}, \underline{\rho}, d_\infty$,
 such that
 \be
 \|\rho(t)-\rho_\infty\|_{H^{-1}}+\|v(t)\|_{\mathbf{H}^1}+\|d(t)-d_\infty\|_{\mathbf{H}^2}\leq C(1+t)^{-\kappa}, \quad \forall\, t \geq
 0, \label{rate}
 \ee
 where $\kappa=\min\left\{\frac{\theta}{1-2\theta}, \frac{\xi}{2}\right\}$ with $\theta \in \left(0,\frac12\right)$ being a constant depending on $f, d_\infty$.
 \et

The question on the uniqueness of asymptotic limit for liquid crystal system \eqref{a1}--\eqref{a3} was raised in \cite{LL95}.
A positive answer was given in \cite{W10} such that for any global classical solution $(v,d)$ of the homogeneous system, the velocity field will decay to zero and the director vector will converge to a steady state that is a solution to the stationary problem \eqref{staa} (cf. \cite{PRS,GW12} for some generalizations). Decay of the velocity field can be obtained by exploring the dissipative nature of the problem and by proper energy estimates, while convergence of the director vector is a nontrivial problem, because the structure of the equilibria set (namely, the set of solutions to \eqref{staa}) can be
quite complicated. The proof in \cite{W10} relies on the so-called \L ojasiewicz--Simon approach \cite{S83}, which turns out to be a useful method
to study convergence of global solutions to equilibria for nonlinear evolution equations (see e.g.,
\cite{Ben, FS, HT01, J981, LD} and the references cited therein). One advantage of the approach is that, one can obtain the convergence
result without studying the structure of the set of equilibria, which is usually difficult when the spacial dimension is
larger than one. 

The nonhomogeneous problem \eqref{0}--\eqref{5} under consideration is much more involved than the homogeneous case. Main difficulties come from those nonlinear couplings between the three equations for density, velocity and director in terms of convection, extra stress term as well as the external force.
Under both assumptions (F1) and (F2), we are able to derive certain (dissipative) basic energy inequalities for problem \eqref{0}--\eqref{5} and furthermore, some specific higher-order differential inequalities in the sprit of \cite{LL95}, which not only provide uniform-in-time estimates for the global strong solutions but also yield decay property of the velocity. Our results on the decay of velocity fields under both types of external forces imply that the dissipations from the viscosity and the relaxation effect in \eqref{3} are strong enough to compensate the effects of external forces and the density fluctuation as well as the interactions between the fluid and the liquid crystal molecules such that the flow will slow down as time goes to infinity. This extends the result on density-dependent incompressible Navier--Stokes equations driven by
a time-independent external force (cf. (F1)) on bounded domains in 2D (cf. \cite{ZK}). For the asymptotically autonomous external force (cf. (F2)), we are able to apply the \L ojasiewicz--Simon approach to prove convergence of the director vector and thus generalize the previous results in \cite{LL95, W10} for the homogeneous liquid crystal system. In this case, we can also obtain $L^1$-integrability of the velocity field, which together with the transport equation \eqref{0} yields convergence of the density function. Besides, by the \L ojasiewicz--Simon approach we can derive some explicit decay rates of the density, velocity field and director vector. We remark that in the case of time-independent external force, the \L ojasiewicz--Simon method seems fail to apply, thus we are only able to show certain sequential convergence of the density function and the director vector. Our results still hold in three-dimensional case provided that bounded global strong solutions of problem \eqref{0}--\eqref{5} can be obtained (this could be verified, for instance, if the initial data and the external forces are sufficiently small).

 The remaining part of this paper is organized as follows. In Section 2,
 we introduce the functional settings, some preliminary results
   as well as some technical lemmas. In Section 3, we derive some dissipative energy inequalities and some specific higher-order differential inequalities for both types of external forces, which enable us to obtain uniform \textit{a priori} estimates and conclude the existence and uniqueness of global strong solutions to problem \eqref{0}{--\eqref{5} (cf. Theorem \ref{main2d}). Section 4 is devoted to the proof of our main results on the long-time behavior of global strong solutions (cf. Theorems \ref{conv1}, \ref{conv}).

\section{Preliminaries}
 As
usual, $L^p(\Omega)$ and $W^{k,p}(\Omega)$ stand for the Lebesgue
and the Sobolev spaces of real valued functions, with the convention
that $H^k(\Omega)= W^{k,2}(\Omega)$. The spaces of vector-valued
functions are denoted by bold letters, correspondingly. Without any
further specification, $\|\cdot\|$ stands for the norm in
$L^2(\Omega)$ or $\mathbf{L}^2(\Omega)$. We shall denote by $C$ the genetic
constants depending on $\lambda, \gamma, \nu, \eta, \Omega$ and the initial
data. Special dependence will be pointed out explicitly in the text
if necessary.

For $1<q<+\infty$, $X^q$ denotes the space that is a completion in $\mathbf{L}^q$ of the set of solenoidal vector fields with coefficients that belong to $\mathcal{V}=\mathbf{C}_0^\infty(\Omega)\cap \{v: \nabla \cdot v=0\}$:
$$ X^q=\{v\in \mathbf{L}^q(\Omega):\ \nabla \cdot v=0\ \text{in}\ \Omega,\  v\cdot \mathbf{n}=0\ \text{on}\ \Gamma\},$$
where $\mathbf{n}$ is the unit outer normal to the boundary (cf. \cite{Da}). In particular, we denote
 \be
 H=X^2=\ \text{the closure of}\ \mathcal{V} \text{\ in\ } \mathbf{L}^2(\Omega),\quad  V=\ \text{the closure of\ } \mathcal{V} \text{\ in\ } \ \mathbf{H}^1_0(\Omega).\non
 \ee
It is well-known that any vector-field with coefficients in $\mathbf{L}^q$ has a Helmholtz decomposition (cf. \cite{SS}). Denote by ${\rm P}_q:\mathbf{L}^q(\Omega)\to X^q$ the projector from $\mathbf{L}^q$ to $X^q$, which is a bounded operator. We can define the Stokes operator $A_q={\rm P}_q(-\Delta)$ with domain $D(A_q)=\mathbf{W}^{2,q}\cap \mathbf{W}^{1,q}_0\cap X^q$.
 The following result holds (cf. e.g., \cite{Da}):
  \bl \label{S}
  Let $\Omega$ be a bounded domain of $\mathbb{R}^2$ with smooth boundary and $d(\Omega)$ be the diameter of $\Omega$. Then following results hold true:

  (i) For any $v\in  \mathbf{W}^{2,q}\cap \mathbf{W}^{1,q}_0$, there is a constant $C=C(q, d(\Omega))$ such that
   $$ \|v\|_{\mathbf{W}^{2,q}}:=\|\nabla ^2 v\|_{\mathbf{L}^q}+d(\Omega)^{-1}\|\nabla v\|_{\mathbf{L}^q}+d(\Omega)^{-2}\|v\|_{\mathbf{L}^q}\leq C\|\nabla^2 v\|_{\mathbf{L}^q}
   $$

   (ii) For $1<q<+\infty$, $f\in \mathbf{L}^q$, the Stokes problem
   $$ -\Delta v+\nabla P=f, \ \text{in}\ \ \Omega, \quad v|_\Gamma=0, $$
   has a unique solution $(v, P)$ in $D(A_q)\times \mathbf{W}^{1,q}$. There exists a constant $C=C(q, d(\Omega))$ such that
   $$
   \|\nabla ^2 v\|_{\mathbf{L}^q}+\|\nabla P\|_{\mathbf{L}^q} \leq C\|f\|_{\mathbf{L}^q}.
   $$
    \el

  Concerning the transport equation for the density function, we have the following result that can be found e.g., in \cite{Da}
  \bl\label{trans}
  Let $v\in L^1(0,T; {\rm Lip})$ be a solenoidal vector field such that $v\cdot \mathbf{n}=0$ on $\Gamma$. For any $\rho_0\in W^{1,q}$ with $q\in [1,+\infty]$, the equation
  $$\rho_t+v\cdot \nabla \rho=0, \quad \rho|_{t=0}=\rho_0(x),$$
  admits a unique solution $\rho\in L^\infty(0, T; W^{1,\infty})\,\cap\, C([0,T]; \cap_{r<\infty}W^{1,r})$ if $q=+\infty$ and $\rho\in C([0,T]; W^{1,q})$ if $1\leq q<+\infty$.
  Besides, the following estimate holds
  \be
   \|\rho(t)\|_{W^{1,q}}\leq e^{\int_0^t\|\nabla v(\tau)\|_{\mathbf{L}^\infty}d\tau}\|\rho_0\|_{W^{1,q}}, \quad \forall\, t\in [0,T].\label{w1qrho}
  \ee
If $\rho_0\in L^p$ for some $p\in [1,+\infty]$, then $\|\rho(t)\|_{L^p}=\|\rho_0\|_{L^p}$ for $t\in [0,T]$.
  \el

An essential characteristic of the director equation (for any given given velocity) is the following weak maximum principle (cf. \cite{LL95}):
\begin{lemma} \label{md}
Suppose $v \in L^\infty(0, T; H) \cap L^2(0, T;
V)$, $d_0\in \mathbf{H}^1(\Omega)$ with $d_0\in \mathbf{H}^\frac32(\Gamma)$ and $|d_0|\leq 1$.
If $d\in L^\infty(0,T; \mathbf{H}^1)\cap
L^2(0, T; \mathbf{H}^2)$ is the weak solution of the initial boundary value
problem
 \bea
&& d_t+ v\cdot \nabla d=\gamma (\Delta d- f(d)), \quad \text{a.e. in\ } \Omega,\non\\
&& d|_\Gamma=d_0(x), \quad (x, t)\in \Gamma\times (0,T)\non\\
&& d|_{t=0}=d_0(x),\non
 \eea
then  $|d(x,t)| \leq 1$, a.e. in $\Omega\times (0, T)$.
\end{lemma}

  Finally, we report some inequalities that will be frequently used in the subsequent proof.

  \bl \label{embed}
  Let $\Omega$ be a bounded domain of $\mathbb{R}^2$ with smooth boundary.
\\  (i) Ladyzhenskaya inequality. $\|v\|_{\mathbf{L}^4}^2\leq C\|\nabla v\|\|v\|, \quad \forall\,v\in \mathbf{H}^1_0(\Omega). $
  \\
  (ii) Agmon inequality. $\|f\|_{L^\infty}^2\leq C\|f\|_{H^2}\|f\|, \quad \forall\,f\in H^2(\Omega). $\\
  (iii) Poincar\'e inequality.  For $1<q<+\infty$, $ \|v\|_{\mathbf{L}^{q}}\leq C\|\nabla v\|_{\mathbf{L}^q}, \quad \forall\,v\in \mathbf{W}^{1,q}_0(\Omega)$.\\
  (iv) Sobolev embeddings. For $1<q<+\infty$, the embedding $H^1\hookrightarrow L^q$ is compact. Besides, $\|f\|_{L^\infty}\leq C\|f\|_{W^{1,q}}, \quad \forall\,f\in W^{1,q}(\Omega)$.\\
  The constant $C$ in the above inequalities may depend on $\Omega$ and $q$.
  \el


\section{ Global strong solutions}
\setcounter{equation}{0}
In order to prove the existence of strong solutions, we can first construct a sequence of approximate solutions $(\rho_m, v_m, d_m)$ within a semi-Galerkin scheme as in \cite[Section 5]{DQS} (we also refer to \cite{LL95} for the case of homogeneous liquid crystal flow and \cite{AKM} for the density-dependent Navier--Stokes equation with external forces).
  To prove the convergence of the approximate solutions, we only need to derive some \textit{a priori} estimates for them. Due to the positivity condition \eqref{bd}, the Galerkin approximation in \cite{DQS} does not rely on the introduction of a viscosity term in the transport equation \eqref{0} like in \cite{LZ09, Tan} and thus it can be used to establish (higher-order) Ladyzhenskaya type energy estimates. Since the calculations for the approximate solutions are (formally) identical to that as we work with smooth solutions, in what follows, we simply perform calculations for smooth solutions to problem \eqref{0}--\eqref{5}.

 \subsection{Dissipative basic energy inequalities and lower-order energy estimates}
 The total energy of problem \eqref{0}--\eqref{5} is defined as follows:
 \be
 \mathcal{E}(t)=\frac{1}{2}\int_\Omega \rho(t)  |v(t)|^2 dx +\frac{\lambda}{2}\|\nabla
 d(t)\|^2+\lambda\int_\Omega F(d(t))dx.\label{Ly}
 \ee
  In analogy to the constant density case (cf. \cite{LL95}) or the nonhomogeneous system without external force (cf. \cite{LZ09, Tan}), our system \eqref{0}--\eqref{5} still has the following \textit{basic
energy inequalities}, which reflect the energy dissipation of the liquid crystal flow .
 \bl[Basic energy inequalities]
Let $(\rho, v,d)$ be a smooth solution of problem \eqref{0}--\eqref{5} on $\Omega\times [0,T]=Q_T$ $(0\leq T\leq +\infty)$.

(i) If $\mathbf{g}$ satisfies (F1), it holds
 \be
 \frac{d}{dt}\tilde{\mathcal{E}}(t)+\nu\|\nabla v(t)\|^2+\lambda\gamma\|\Delta
 d(t)-f(d(t))\|^2=0,
 \quad 0\leq t\leq T,\label{lya}
 \ee
 where $$ \tilde{\mathcal{E}}(t)=\mathcal{E}(t)-\int_\Omega \rho\phi dx.$$

 (ii) If $\mathbf{g}$ satisfies (F2), it holds
 \be
 \frac{d}{dt}\mathcal{E}(t)+\frac{\nu}{2}\|\nabla v(t)\|^2+\lambda\gamma\|\Delta
 d(t)-f(d(t))\|^2\leq \frac{C_P^2\bar{\rho}^2}{2\nu}\|\mathbf{g}\|^2,
 \quad 0\leq t\leq T.\label{lya1}
 \ee
 \el
 \begin{proof}
 Multiplying \eqref{0}, \eqref{1} and \eqref{3} by $\frac12|v|^2$, $v$ and $\lambda(-\Delta d + f(d))$, respectively, integrating
over $\Omega$, we infer from the boundary conditions \eqref{5}  and integration by parts that
 \be
  \frac12 \int_\Omega \rho_t |v|^2 dx  = -\frac12 \int_\Omega \nabla \cdot(\rho v)|v|^2 dx= \frac12 \int_\Omega \rho v\cdot \nabla |v|^2dx, \label{b0}
  \ee
 \bea && \frac12\int_\Omega  \rho \frac{d}{dt}|v|^2 dx+\nu \|\nabla v\|^2\non\\
 &&\ \ =-\frac12 \int_\Omega \rho v\cdot \nabla |v|^2dx +\int_\Omega \rho \mathbf{g} \cdot v dx -\lambda \int (\Delta d\cdot\nabla d )\cdot v dx,\label{b1}
 \eea
 \be
 \frac{d}{dt}\left(\frac{\lambda}{2}\|\nabla d\|^2+\lambda\int_\Omega F(d)dx\right)+\lambda\gamma\|-\Delta f+f(d)\|^2+\lambda \int_\Omega (v\cdot \nabla d)\cdot \Delta d dx=0,\label{b2}
 \ee
 where we have used the facts (cf. \cite{LL95})
 \bea
 && \nabla \cdot (\nabla d\odot \nabla d)=\frac12 \nabla |\nabla d|^2+  \Delta d\cdot \nabla d,\label{cann1}\\
 && \int_\Omega \nabla P \cdot vdx=\int_\Omega \nabla |\nabla d|^2 \cdot v dx=\int_\Omega v\cdot \nabla F(d) dx=0.\label{cann2}
 \eea
 Adding \eqref{b0}--\eqref{b2} together, we can see that
 \be \frac{d}{dt}\mathcal{E}(t)+\nu\|\nabla v(t)\|^2+\lambda\gamma\|\Delta
 d(t)-f(d(t))\|^2=\int_\Omega \rho \mathbf{g} \cdot v dx. \label{bb}
 \ee
 If $\mathbf{g}$ satisfies (F1), using the idea in \cite{ZK}, we multiply the transport equation \eqref{0} by $-\phi$ and integrate over $\Omega$ to get
 \be
 -\int_\Omega  \rho_t \phi dx  = \int_\Omega \nabla \cdot(\rho v)\phi dx= -\int_\Omega \rho \nabla\phi \cdot v dx.\label{b3}
 \ee
Adding \eqref{bb} with \eqref{b3} and noticing that $\phi$ is independent of time, we arrive at our conclusion \eqref{lya}. On the other hand, If $\mathbf{g}$ satisfies (F2), using the H\"older inequality and Poincar\'e inequality, we infer that
 \be \left|\int_\Omega \rho \mathbf{g} \cdot v dx \right|\leq \|\rho\|_{L^\infty}\|\mathbf{g}\|\|v\|\leq C_P\bar{\rho}\|\mathbf{g}\|\|\nabla v\|\leq \frac{\nu}{2}\|\nabla v\|^2+ \frac{C_P^2\bar{\rho}^2}{2\nu}\|\mathbf{g}\|^2,
 \ee
 which together with \eqref{bb} yields \eqref{lya1}.
 \end{proof}

  \bp \label{v0h1}
  Under the assumptions of Theorem \ref{main2d}, the following estimates hold
  \bea
  && 0<\underline{\rho}\leq \rho(x,t)\leq \bar{\rho},\quad \forall\,
 t\geq 0,\label{bdr}\\
 && |d(x,t)|\leq 1, \quad \forall\, t\geq 0,\label{di1}\\
  && \|v(t)\|+\|d(t)\|_{\mathbf{H}^1}\leq C,\quad  \forall\, t\geq 0,\label{bd1}\\
  && \int_0^{+\infty} (\nu\|\nabla v(t)\|^2+\lambda\gamma\|\Delta
 d(t)-f(d(t))\|^2) dt\leq C,\label{intA}
  \eea
  where $C$ is a constant depending  on $\|v_0\|$, $\|d_0\|_{\mathbf{H}^1}$, $\eta$, $\underline{\rho}$, $\bar{\rho}$, $\Omega$ and also $\|\phi\|_{H^1}$ (under (F1)) or $\|\mathbf{g}\|_{L^2(0,+\infty;\, \mathbf{L}^2)}$ (under (F2)).
  \ep
  \begin{proof}
  \eqref{bdr} easily follows from the characteristics method \cite{La}, while \eqref{di1} is a consequence of the weak maximum principle for the director equation (see Lemma \ref{md}).   If $\mathbf{g}$ satisfies (F1), since
 \be
 \left|\int_\Omega \rho \phi dx\right|\leq |\Omega|\bar{\rho}\|\phi\|_{L^1}<+\infty,\label{rp}
 \ee
 we can see that $\tilde{\mathcal{E}}(t)\geq -|\Omega|\bar{\rho}\|\phi\|_{L^1}$ for all $t\geq 0$.  The required uniform estimates follow from \eqref{lya} and \eqref{rp}. If $\mathbf{g}$ satisfies (F2), by integrating \eqref{lya1} with respect to time, we arrive at the conclusion.
   \end{proof}

\subsection{Higher-order energy estimates}
Denote
 \be
 A(t)=\nu \|\nabla v(t)\|^2 + \|\Delta d(t)-f(d(t))\|^2.\label{A}
 \ee
 \bl \label{he2d} The following inequality holds for smooth solutions $(\rho, v,d)$ to problem \eqref{0}--\eqref{5}:
 \bea && \frac{d}{dt}A(t)+\|\rho^\frac12 v_t(t)\|^2+\gamma\|\nabla(\Delta d(t)-f(d(t)))\|^2\non\\
 & \leq& C(A^2(t)+A(t))+C\|\mathbf{g}\|^2, \ \  \forall
 \,
 t> 0,\label{he}
 \eea
 where $C$ is a constant depending on $\|v_0\|, \|d_0 \|_{\mathbf{H}^1}, \eta, \underline{\rho}, \bar{\rho}, \nu, \Omega$ and also $\|\phi\|_{H^1}$ (under (F1)) or $\|\mathbf{g}\|_{L^2(0,+\infty;\, \mathbf{L}^2)}$ (under (F2)).
 \el
\begin{proof}
Using equations  \eqref{0}--\eqref{3} and the facts \eqref{cann1}, \eqref{cann2}, we compute that
\bea
 && \frac{\nu}{2}\frac{d}{dt}\|\nabla v\|^2 = -\int_\Omega \nu \Delta v\cdot v_t dx  \non\\
&=& -\int_\Omega \rho |v_t|^2dx-\int_\Omega \rho(v\cdot \nabla v)\cdot v_tdx\non\\
&& -\lambda \int_\Omega [ (\Delta d-f(d))\cdot \nabla d]\cdot v_tdx  +\int_\Omega \rho \mathbf{g} \cdot v_t dx
\label{A1}
\eea
and
\bea
&& \frac{1}{2}\frac{d}{dt} \|\Delta d-f(d)\|^2\non\\
&=&  \int_\Omega (\Delta d_t-f'(d)d_t)\cdot (\Delta d-f(d)) dx\non\\
&=& - \int \Delta (v\cdot \nabla d) \cdot (\Delta d-f(d)) dx -\gamma \|\nabla (\Delta d-f(d))\|^2\non\\
&& - \int_\Omega f'(d)[-v\cdot \nabla d+\gamma(\Delta d-f(d))](\Delta d-f(d)) dx\non\\
&=& - \int (\Delta v \cdot \nabla d) (\Delta d-f(d)) dx-\int_\Omega [(v\cdot \nabla)\Delta d +2\nabla v\nabla^2 d] \cdot (\Delta d-f(d)) dx\non\\
&&  - \gamma \|\nabla (\Delta d-f(d))\|^2 -\gamma \int_\Omega f'(d)(\Delta d-f(d))\cdot (\Delta d-f(d)) dx.\label{A2}
\eea
Adding \eqref{A1} and \eqref{A2} together, we have
\bea
&&\frac12\frac{d}{dt}\left(\nu \|\nabla v\|^2 + \|\Delta d-f(d)\|^2\right)+\int_\Omega \rho |v_t|^2dx+\gamma \|\nabla (\Delta d-f(d))\|^2\non\\
&=& -\int_\Omega \rho(v\cdot \nabla v)\cdot v_tdx- \int_\Omega [ (\Delta d-f(d))\cdot\nabla d]\cdot (\Delta v+ \lambda v_t)dx\non\\
&& +\int_\Omega \rho \mathbf{g} \cdot v_t dx-\int_\Omega (v\cdot \nabla)\Delta d \cdot (\Delta d-f(d)) dx\non\\
&& -2\int_\Omega (\nabla v\nabla^2 d)\cdot (\Delta d-f(d)) dx-\gamma \int_\Omega f'(d)(\Delta d-f(d))\cdot (\Delta d-f(d)) dx\non\\
&:=& \sum_{m=1}^6 I_m.\label{dA}
\eea
Next, we estimate $I_1, ..., I_6$ term by term.
Using the fact \eqref{cann1}, we rewrite \eqref{1} as
 \be
 -\nu \Delta v+\nabla \left(P+\frac{\lambda}{2} |\nabla d|^2+\lambda F(d)\right)=-\rho(v_t+v\cdot\nabla v)-\lambda
 (\Delta d-f(d))\cdot \nabla d+\rho\mathbf{g}.\label{vvv}
 \ee
By Lemma \ref{S}, we get
 \bea
 \|v\|_{\mathbf{H}^2}
 &\leq& C(\|\rho v_t\|+\|\rho(v\cdot\nabla)v\|+\| (\Delta d-f(d))\cdot \nabla d\|+\|\rho\mathbf{g}\|)\non\\
 &\leq& C(\bar{\rho}^\frac12\|\rho^\frac12 v_t\|+\bar{\rho}\|v\|_{\mathbf{L}^4}\|\nabla v\|_{\mathbf{L}^4}+\bar{\rho}\|\mathbf{g}\|)+ C\|(\Delta d-f(d))\cdot \nabla d\|\non\\
 &\leq& C(\|\rho^\frac12v_t\|+C\|v\|^\frac12\|\nabla v\|\|v\|_{\mathbf{H}^2}^\frac12+\|\mathbf{g}\|)+ C\|(\Delta d-f(d))\cdot \nabla d\|\non\\
 &\leq& \frac12\|v\|_{\mathbf{H}^2}+ C(\|\rho^\frac12 v_t\|+\|\mathbf{g}\|+\|(\Delta d-f(d))\cdot \nabla d\|)+ C\|\nabla v\|^2,\label{vh2e}
 \eea
 namely,
  \be
 \|v\|_{\mathbf{H}^2}
 \leq C(\|\rho^\frac12 v_t\|+\|\mathbf{g}\|+\|(\Delta d-f(d))\cdot \nabla d\|)+ C\|\nabla v\|^2.\label{vh2}
 \ee
 On the other hand, we have
 \bea
 && \|(\Delta d-f(d))\cdot \nabla d\|\non\\
 &\leq& \|\nabla d\|_{\mathbf{L}^4}\|\Delta d-f(d)\|_{\mathbf{L}^4}\non\\
 &\leq& C(\|\Delta d\|^\frac12\|\nabla d\|^\frac12+\|\nabla d\|) \|\Delta d-f(d)\|^\frac12\|\nabla (\Delta d-f(d))\|^\frac12\non\\
 &\leq& C(\|\Delta d-f(d)\|^\frac12+1)\|\Delta d-f(d)\|^\frac12\|\nabla (\Delta d-f(d))\|^\frac12.\label{dd}
 \eea
 Then we infer from \eqref{vh2}, \eqref{dd} and the Young inequality that
 \bea
 I_1&\leq& \frac{1}{16}\int_\Omega \rho|v_t|^2dx+4\int_\Omega \rho|v\cdot \nabla v|^2dx\non\\
 &\leq& \frac{1}{16}\int_\Omega \rho|v_t|^2dx+4\bar{\rho} \|v\|_{\mathbf{L}^4}^2\|\nabla v\|_{\mathbf{L}^4}^2\non\\
 &\leq& \frac{1}{16}\int_\Omega \rho|v_t|^2dx+C \|v\|\|\nabla v\|(\|\nabla v\|\|\Delta v\|+\|\nabla v\|^2)\non\\
 &\leq& \frac{1}{16}\int_\Omega \rho|v_t|^2dx+C\|v\|_{\mathbf{H}^2}\|\nabla v\|^2+ C\|\nabla v\|^3\non\\
 &\leq& \frac{1}{16}\int_\Omega \rho|v_t|^2dx+ C\|\rho^\frac12 v_t\|\|\nabla v\|^2+C(\|\mathbf{g}\|+\|\nabla v\|^2)\|\nabla v\|^2+C\|\nabla v\|^3\non\\
 && +C(\|\Delta d-f(d)\|^\frac12+1)\|\Delta d-f(d)\|^\frac12\|\nabla (\Delta d-f(d))\|^\frac12\|\nabla v\|^2\non\\
 &\leq& \frac18\int_\Omega \rho|v_t|^2dx+\frac{\gamma}{8} \|\nabla (\Delta d-f(d))\|^2+C\|\nabla v\|^2+C\|\nabla v\|^4\non\\
 && +C\|\Delta d-f(d)\|^2+C\|\Delta d-f(d)\|^4+C\|\mathbf{g}\|^2,\non
 \eea
  \bea
 I_2&\leq& C(\|\Delta v\|+\underline{\rho}^\frac12 \|\rho^\frac12 v_t\|)\|(\Delta d-f(d))\cdot \nabla d\|\non\\
 &\leq&  C(\|\rho^\frac12 v_t\|+\|\mathbf{g}\|+\|(\Delta d-f(d))\cdot \nabla d\|+\|\nabla v\|^2)\|(\Delta d-f(d))\cdot \nabla d\|\non\\
 &\leq& C(\|\rho^\frac12v_t\|+\|\mathbf{g}\|+\|\nabla v\|^2)\non\\
 &&\ \times (\|\Delta d-f(d)\|^\frac12+1)\|\Delta d-f(d)\|^\frac12\|\nabla (\Delta d-f(d))\|^\frac12\non\\
 &&+C (\|\Delta d-f(d)\|^2+\|\Delta d-f(d)\|)\|\nabla (\Delta d-f(d))\|\non\\
 &\leq& \frac18\int_\Omega \rho|v_t|^2dx+ \frac{\gamma}{8} \|\nabla (\Delta d-f(d))\|^2 \non\\
 && +C\|\Delta d-f(d)\|^4+C\|\Delta d-f(d)\|^2+C\|\nabla v\|^4+C\|\mathbf{g}\|^2.\non
 \eea
 Concerning $I_3$, it follows that
 \be
 I_3\leq \bar{\rho}^\frac12\|\mathbf{g}\|\|\rho^\frac12 v_t\|\leq \frac{1}{8}\|\rho v_t\|^2+ 2\bar{\rho}\|\mathbf{g}\|^2.\non
 \ee
  Next, using \eqref{bd1}, we get
 \be
 \|f'(d)\nabla d\|\leq C(\|d\|_{\mathbf{L}^8}^2+1)\|\nabla d\|_{\mathbf{L}^4}\leq C(1+\|\Delta d-f(d)\|^\frac12),\non
 \ee
 which implies that
 \bea
 I_4&\leq& C\|v\|_{\mathbf{L}^4}\|\nabla \Delta d\|\|\Delta d-f(d)\|_{\mathbf{L}^4}\non\\
 &\leq& C\|v\|^\frac12\|\nabla v\|^\frac12(\|\nabla (\Delta d-f(d))\|+\|f'(d)\nabla d\|) \|\Delta d-f(d)\|^\frac12\|\nabla (\Delta d-f(d))\|^\frac12\non\\
 &\leq& C\|v\|^\frac12\|\nabla v\|^\frac12\|\Delta d-f(d)\|^\frac12\|\nabla (\Delta d-f(d))\|^\frac32\non\\
 && +C\|v\|^\frac12\|\nabla v\|^\frac12(\|\Delta d-f(d)\|^\frac12+\|\Delta d-f(d)\|)\|\nabla (\Delta d-f(d))\|^\frac12\non\\
 &\leq& \frac{\gamma}{8}\|\nabla (\Delta d-f(d))\|^2+C\|\Delta d-f(d)\|^4+C\|\Delta d-f(d)\|^2\non\\
 && +C\|\nabla v\|^4+C\|\nabla v\|^2.\non
 \eea
 From the elliptic estimate and \eqref{bd1}, we have
 \bea
 \|d\|_{\mathbf{H}^2}&\leq& C\left(\|\Delta d\|+\|d\|_{\mathbf{H}^\frac32(\Gamma)}\right)\leq C(\|\Delta d-f(d)\|+\|f(d)\|+\|d_0\|_{\mathbf{H}^2})\non\\
 &\leq&
 C(\|\Delta d-f(d)\|+1).\label{dh2}
 \eea
 Then it follows from \eqref{vh2}, \eqref{dd}, \eqref{dh2} and Young inequality that
 \bea
 I_5&\leq& C\|\nabla v\|_{\mathbf{L}^4}\|d\|_{\mathbf{H}^2}\|\Delta d-f(d)\|_{\mathbf{L}^4}\non\\
 &\leq& C\|\nabla v\|^\frac12\|v\|_{\mathbf{H}^2}^\frac12(\|\Delta d-f(d)\|+1)\|\Delta d-f(d)\|^\frac12\|\nabla (\Delta d-f(d))\|^\frac12\non\\
 &\leq& C\|\nabla v\|^\frac12(\|\rho^\frac12 v_t\|^\frac12 +\|\mathbf{g}\|^\frac12 +\|\nabla v\|)(\|\Delta d-f(d)\|+1)\non\\
 && \times\|\Delta d-f(d)\|^\frac12\|\nabla (\Delta d-f(d))\|^\frac12\non\\
 && + C\|\nabla v\|^\frac12(\|\Delta d-f(d)\|^\frac54+1)\|\Delta d-f(d)\|^\frac34\|\nabla (\Delta d-f(d))\|^\frac34\non \\
 &\leq & \frac18\int_\Omega \rho|v_t|^2dx+\frac{\gamma}{8}\|\nabla (\Delta d-f(d))\|^2+C\|\mathbf{g}\|^2\non\\
 &&+ C\|\nabla v\|^4+C\|\nabla v\|^2+C\|\Delta d-f(d)\|^4+C\|\Delta d-f(d)\|^2.\non
 \eea
 Finally, for $I_6$, we have
 \bea
 I_6&\leq& C\|f'(d)\|_{\mathbf{L}^4}\|\Delta d-f(d)\|_{\mathbf{L}^4}\|\Delta d-f(d)\|\non\\
 &\leq& C(\|d\|_{\mathbf{L}^8}^2+1)\|\nabla (\Delta d-f(d))\|^\frac12\|\Delta d-f(d)\|^\frac32\non\\
 &\leq& \frac{\gamma}{8}\|\nabla (\Delta d-f(d))\|^2+ C\|\Delta d-f(d)\|^2.\non
 \eea
 Collecting the estimates for $I_1,..., I_6$, we infer from \eqref{dA} that
 \bea
 \frac12\frac{d}{dt}A(t)&\leq&   -\frac12  \int_\Omega \rho|v_t|^2dx-\frac{\gamma}{2} \|\nabla (\Delta d-f(d))\|^2 + C\|\nabla v\|^4+C\|\nabla v\|^2\non\\
 &&+C\|\Delta d-f(d)\|^4+C\|\Delta d-f(d)\|^2+C\|\mathbf{g}\|^2,\non
 \eea
 which yields our conclusion \eqref{he}.
\end{proof}
\bp \label{es12}
 Under the assumptions of Theorem \ref{main2d}, the following uniform estimates hold for any $t\geq 0$:
 \bea
 && \|v(t)\|_{\mathbf{H}^1}+ \|d(t)\|_{\mathbf{H}^2}\leq C, \label{ubdd}\\
 && \sup_{t\geq 0} \int_t^{t+1} \|\rho^\frac12 v_t(\tau)\|^2+\|\nabla(\Delta d(\tau)-f(d(\tau)))\|^2 d\tau \leq C, \label{intbd1}
 \eea
 where $C$ is a constant depending on $\nu, \|v_0\|_{\mathbf{H}^1}, \|d_0
 \|_{\mathbf{H}^2}, \eta, \underline{\rho}, \bar{\rho}$ and also $\|\phi\|_{H^1}$ (under (F1)) or $\|\mathbf{g}\|_{L^2(0,+\infty;\, \mathbf{L}^2)}$ (under (F2)).
  \ep
 \begin{proof}
 For both the cases (F1) and (F2), \eqref{intA} implies that
 $$\int_0^{+\infty} A(t) dt<+\infty.$$
 The uniform
 bound \eqref{ubdd} follows from Lemma \ref{v0h1}, the higher-order energy inequality \eqref{he} and the uniform Gronwall lemma \cite[Lemma III.1.1]{Te}.
  Integrating \eqref{he} from $t$ to $t+1$, we can conclude \eqref{intbd1} from \eqref{ubdd}.
 \end{proof}

 Denote the quantity
 \be
 B(t)=\int_\Omega \rho|v_t(t)|^2dx+\|\nabla d_t(t)\|^2.\label{B}
 \ee
 We can derive the following higher-order differential inequality:
 \bl
 If $\mathbf{g}$ satisfies (F1), then following inequality holds for smooth solutions $(\rho, v,d)$ to problem \eqref{0}--\eqref{5}:
 \be
 \frac{d}{dt}B(t)+\nu \|\nabla v_t\|^2+\gamma \|\Delta d_t\|^2\leq C(B^2(t)+B(t))+C\|\nabla v\|^2, \ \ \forall\, t\geq 0,\label{eB1}
 \ee
 where $C$ is a constant depending on $\nu, \|v_0\|_{\mathbf{H}^1}, \|d_0
 \|_{\mathbf{H}^2}, \eta, \underline{\rho}, \bar{\rho}$ and $\|\phi\|_{H^2}$.
 On the other hand, if $\mathbf{g}$ satisfies (F2), then we have:
 \be
 \frac{d}{dt}B(t)+\nu \|\nabla v_t\|^2+\gamma \|\Delta d_t\|^2\leq C(B^2(t)+B(t))+C(\|\mathbf{g}\|_{\mathbf{H}^1}^2+\|\mathbf{g}_t\|^2+\|\nabla v\|^2), \ \ \forall\, t\geq 0,\label{eB2}
 \ee
 where $C$ is a constant depending on $\nu, \|v_0\|_{\mathbf{H}^1}, \|d_0
 \|_{\mathbf{H}^2}, \eta,  \underline{\rho}, \bar{\rho}$ and $\|\mathbf{g}\|_{L^2(0,+\infty;\, \mathbf{L}^2)}$.
 \el
 \begin{proof}
 Taking temporal derivative of \eqref{1} and taking $L^2$ inner product of the resultant with $v_t$, then using the fact $\nabla \cdot v=\nabla \cdot v_t=0$ and the transport equation \eqref{0}, after integration by parts (keeping in mind that $v|_\Gamma=v_t|_\Gamma=0$), we get
 \bea
 && \frac12\frac{d}{dt} \int_\Omega \rho|v_t|^2dx + \nu \|\nabla v_t\|^2\non\\
 &=& \frac12\int_\Omega \rho_t |v_t|^2 dx-\frac12 \int_\Omega \rho v\cdot \nabla (|v_t|^2) dx-\int_\Omega \rho (v_t\cdot \nabla) v \cdot v_t dx \non\\
 && -\int_\Omega \rho_t(v_t+v\cdot\nabla v)\cdot v_t dx-\int_\Omega \nabla P_t \cdot v_t dx+\int_\Omega (\rho_t\mathbf{g}+\rho\mathbf{g}_t) \cdot v_t dx\non\\
 && +2\lambda\int_\Omega (\nabla d_t\odot \nabla d):\nabla v_t dx\non\\
 &=& -\int_\Omega \rho_t|v_t|^2dx -\int_\Omega \rho_t(v\cdot\nabla v)\cdot v_t dx-\int_\Omega \rho (v_t\cdot \nabla v) \cdot v_t dx\non\\
   &&+\int_\Omega (\rho_t\mathbf{g}+\rho\mathbf{g}_t) \cdot v_t dx+2\lambda\int_\Omega (\nabla d_t\odot \nabla d):\nabla v_t dx\non\\
   &=& -\int_\Omega \rho v\cdot \nabla (|v_t|^2)dx -\int_\Omega \rho v \cdot \nabla [ (v\cdot\nabla v)\cdot v_t] dx -\int_\Omega \rho (v_t\cdot \nabla v) \cdot v_t dx
   \non\\
   && +\int_\Omega \rho v\cdot \nabla (\mathbf{g} \cdot v_t) dx+ \int_\Omega \rho\mathbf{g}_t \cdot v_t dx +2\lambda\int_\Omega (\nabla d_t\odot \nabla d):\nabla v_t dx
   \non\\
   &:=&\sum_{m=1}^{6}J_m.\label{dvt}
 \eea
 First, we consider the case that $\mathbf{g}$ satisfies the assumption (F1). In this case, $J_5=0$. Using the uniform estimates in Proposition \ref{es12}, we have
 \bea
 J_1&\leq & C\bar{\rho} \|\nabla v_t\|\|v_t\|_{\mathbf{L}^4}\|v\|_{\mathbf{L}^4} \leq C \|\nabla v_t\|(\|\nabla v_t\|^\frac12\|v_t\|^\frac12+\|v_t\|)\non\\
 &\leq& \frac{\nu}{12}\|\nabla v_t\|^2+C\|v_t\|^2.\non
 \eea
  Then for $J_2$, we infer from \eqref{vh2}, \eqref{dd} and \eqref{ubdd} that
 \bea
 J_2&\leq& \bar{\rho}(\|\nabla v\|_{\mathbf{L}^4}^2\|v\|_{\mathbf{L}^4}\|v_t\|_{\mathbf{L}^4}+\|v\|_{\mathbf{L}^8}^2\|v\|_{\mathbf{H}^2}\|v_t\|_{\mathbf{L}^4}+\|v\|_{\mathbf{L}^8}^2\|\nabla v\|_{\mathbf{L}^4}\|\nabla v_t\|)\non\\
 &\leq & C[ (\|\Delta v\|\|\nabla v\|+\|\nabla v\|^2)\|\nabla v\|+\|v\|_{\mathbf{H}^2}\|\nabla v\|^2\non\\
 && + \|\nabla v\|^2(\|\Delta v\|^\frac12\|\nabla v\|^\frac12+\|\nabla v\|)]\|\nabla v_t\|\non\\
 &\leq& C\|\nabla v\|\|\nabla v_t\|+ C\|v\|_{\mathbf{H}^2}\|\nabla v\|\|\nabla v_t\|\non\\
 &\leq&  C\|\nabla v\|\|\nabla v_t\|+ C(\|\rho^\frac12 v_t\|+\|\mathbf{g}\|+\|\nabla v\|^2)\|\nabla v\|\|\nabla v_t\|\non\\
 &&+C (\|\Delta d-f(d)\|^\frac12+1)\|\Delta d-f(d)\|^\frac12\|\nabla (\Delta d-f(d))\|^\frac12\|\nabla v\|\|\nabla v_t\|\non\\
 &\leq& \frac{\nu}{12}\|\nabla v_t\|^2+C\|\nabla (\Delta d-f(d))\|^2\|\nabla v\|^2+ C\|\rho^\frac12 v_t\|^2\|\nabla v\|^2\non\\
 && + C(1+\|\mathbf{g}\|^2) \|\nabla v\|^2.\label{J2}
 \eea
 When $\mathbf{g}=\nabla \phi$, then the last term on the right-hand side of \eqref{J2} is controlled by $C\|\nabla v\|^2$. Next,
 \bea
 J_3&\leq& \bar{\rho}\|\nabla v\|\|v_t\|_{\mathbf{L}^4}^2\leq C\|\nabla v_t\|\|v_t\|\leq \frac{\nu}{12}\|\nabla v_t\|^2+ C\|v_t\|^2,\non\\
 J_4&=& \int_\Omega \rho v\cdot \nabla (\nabla \phi \cdot v_t) dx\non\\
 &\leq& \bar{\rho}(\|v\|_{\mathbf{L}^4}\|v_t\|_{\mathbf{L}^4}\|\phi\|_{H^2}+\|\nabla v_t\|\|v\|_{\mathbf{L}^4}\|\nabla \phi\|_{\mathbf{L}^4})\non\\
 &\leq& C\|\nabla v\|\|\nabla v_t\|\non\\
 &\leq& \frac{\nu}{12}\|\nabla v_t\|^2+ C\|\nabla v\|^2.\non
  \eea
 Finally, for $J_6$, we have
 \bea
 J_6&\leq& C\|\nabla d\|_{L^\infty}\|\nabla d_t\|\|\nabla v_t\|\non\\
 &\leq& C\|\nabla d\|_{H^2}\|\nabla d\|\|\nabla d_t\|\|\nabla v_t\|\non\\
 &\leq& C(\|\nabla (\Delta d-f(d))\|+1)\|\nabla d_t\|\|\nabla v_t\|\non\\
 &\leq& \frac{\nu}{12}\|\nabla v_t\|^2+ C(\|\nabla (\Delta d-f(d))\|^2+1)\|\nabla d_t\|^2.\non
 \eea
 Taking temporal derivative of \eqref{3} and taking $L^2$ inner product of the resultant with $-\Delta d_t$, we obtain that
 \bea
 && \frac12\frac{d}{dt}\|\nabla d_t\|^2+\gamma \|\Delta d_t\|^2\non\\
 &=& \int_\Omega (v_t\cdot\nabla d)\cdot \Delta d_t dx+\int_\Omega (v\cdot \nabla d_t)\cdot \Delta d_t dx +\gamma \int_\Omega f'(d) d_t\cdot \Delta d_t dx\non\\
 &:=& J_7+J_8+J_9,\label{d3t}
 \eea
 where the right-hand side of \eqref{d3t} can be estimated as follows
 \bea
 J_7&\leq&
  \|v_t\|_{\mathbf{L}^4}\|\nabla d\|_{\mathbf{L}^4}\|\Delta d_t\|\leq C\|\nabla v_t\|^\frac12\|v_t\|^\frac12\|\Delta d_t\|\non\\
  &\leq& \frac{\nu}{12}\|\nabla v_t\|^2+\frac{\gamma}{6}\|\Delta d_t\|+ C\|v_t\|^2,\non
  \\
 J_8&\leq& \|v\|_{\mathbf{L}^4}\|\nabla d_t\|_{\mathbf{L}^4}\|\Delta d_t\|\leq C (\|\nabla d_t\|^\frac12\|\Delta d_t\|^\frac12+\|\nabla d_t\|)\|\Delta d_t\|\non\\
 &\leq& \frac{\gamma}{6}\|\Delta d_t\|^2+ C\|\nabla d_t\|^2,\non\\
 J_9&\leq& \gamma \|f'(d)\|_{\mathbf{L}^\infty}\|d_t\|\|\Delta d_t\|\leq \frac{\gamma}{6}\|\Delta d_t\|^2+ C\|\nabla d_t\|^2.\non
 \eea
  It follows from \eqref{3} that
 \be
 \nabla d_t=-\nabla (v\cdot\nabla d)+\gamma \nabla (\Delta d-f(d)).\label{dtd}
 \ee
 Using the uniform estimate \eqref{ubdd} and Sobolev embedding theorem, we have
 \bea
\|\nabla (v\cdot\nabla d)\| &\leq& C\|\nabla v\|\|\nabla d\|_{\mathbf{L}^\infty}+C\|v\|_{\mathbf{L}^4}\|\nabla^2 d\|_{\mathbf{L}^4}\non\\
 &\leq& C\|\nabla v\|\|\nabla d\|_{\mathbf{H}^2}^\frac12\|\nabla d\|^\frac12+ C\|\nabla v\|(\|\nabla \Delta d\|^\frac12\|\Delta d\|^\frac12+\|\Delta d\|)\non\\
 &\leq& C\|\nabla v\|\|\nabla \Delta d\|^\frac12+C\|\nabla v\|\non\\
 &\leq& C\|\nabla v\|\|\nabla (\Delta d-f(d))\|^\frac12+C\|\nabla v\|,
 \eea
 which together with the equality \eqref{dtd} implies that
 \bea
  \|\nabla (\Delta d-f(d))\| &\leq& \frac{1}{\gamma}\|\nabla d_t\|+\frac{1}{\gamma}\|\nabla (v\cdot\nabla d)\|\non\\
  &\leq& \frac{1}{\gamma}\|\nabla d_t\|+C\|\nabla v\|\|\nabla (\Delta d-f(d))\|^\frac12+C\|\nabla v\|\non\\
  &\leq& \frac{1}{\gamma}\|\nabla d_t\|+\frac12\|\nabla (\Delta d-f(d))\|+C\|\nabla v\|.\label{eqv1}
 \eea
 As a result, we get
 \bea
 \|\nabla d_t\|&\leq& \|\nabla (v\cdot\nabla d)\|+\gamma \|\nabla (\Delta d-f(d))\|\non\\
 &\leq& C\|\nabla v\|\|\nabla (\Delta d-f(d))\|^\frac12+C\|\nabla v\|+\gamma \|\nabla (\Delta d-f(d))\|\non\\
 &\leq& (1+\gamma) \|\nabla (\Delta d-f(d))\|+C\|\nabla v\|.\label{eqv2}
 \eea
 Collecting the above estimates for $J_1, ..., J_9$ and using the relations \eqref{eqv1}, \eqref{eqv2}, we conclude from \eqref{dvt} and \eqref{d3t} that
 \bea
 && \frac{d}{dt}B(t)+\nu \|\nabla v_t\|^2+\gamma \|\Delta d_t\|^2\non\\
 &\leq& C\|\nabla v\|^2(\|\rho^\frac12 v_t\|^2+\|\nabla d_t\|^2)+C\|\nabla d_t\|^4+C(\|v_t\|^2+\|\nabla d_t\|^2)+C\|\nabla v\|^2,\non
 \eea
 which yields our conclusion \eqref{eB1}.

 Now if $\mathbf{g}$ satisfies (F2), we only have to re-estimate the terms $J_2$, $J_4$ and $J_5$ using the uniform estimate \eqref{ubdd}.
 It follows from \eqref{J2} that
 \bea
 J_2&\leq& \frac{\nu}{12}\|\nabla v_t\|^2+C\|\nabla (\Delta d-f(d))\|^2\|\nabla v\|^2+ C\|\rho^\frac12 v_t\|^2\|\nabla v\|^2+ C(1+\|\mathbf{g}\|^2) \|\nabla v\|^2\non\\
 &\leq& \frac{\nu}{12}\|\nabla v_t\|^2+C\|\nabla (\Delta d-f(d))\|^2\|\nabla v\|^2+ C\|\rho^\frac12 v_t\|^2\|\nabla v\|^2+ C\|\mathbf{g}\|^2+C \|\nabla v\|^2.
 \non
 \eea
 Moreover,
 \bea
 J_4 &\leq& \bar{\rho}(\|v\|_{\mathbf{L}^4}\|v_t\|_{\mathbf{L}^4}\|\mathbf{g}\|_{\mathbf{H}^1}+\|\nabla v_t\|\|v\|_{\mathbf{L}^4}\|\mathbf{g}\|_{\mathbf{L}^4})\non\\
 &\leq& C\|\mathbf{g}\|_{\mathbf{H}^1}\|\|\nabla v\|\|\nabla v_t\|\non\\
 &\leq& \frac{\nu}{24}\|\nabla v_t\|^2+ C\|\mathbf{g}\|_{\mathbf{H}^1}^2,\non
 \\
 J_5&\leq&\bar{\rho} \|\mathbf{g}_t\|\|v_t\|\leq \frac{\nu}{24}\|\nabla v_t\|^2+C\|\mathbf{g}_t\|^2.\non
 \eea
 By the above estimates, we easily arrive at \eqref{eB2}.
  \end{proof}
 \bp \label{est11}
  Under the assumptions of Theorem \ref{main2d}, the following uniform estimates hold for any $t\geq 0$,
 \bea
 && \|v_t(t)\|+ \|d_t(t)\|_{\mathbf{H}^1}\leq C, \label{ubdd1}\\
 && \|v(t)\|_{\mathbf{H}^2}+\|d(t)\|_{\mathbf{H}^3}\leq C, \label{ubdd2}\\
 && \sup_{t\geq 0} \int_t^{t+1} (\|\nabla v_t(\tau)\|^2+ \|\Delta d_t(\tau)\|^2) d\tau \leq C, \label{intbd2}\\
 && \sup_{t\geq 0} \int_t^{t+1} \|v(\tau)\|_{\mathbf{W}^{2,q}}^2 d\tau \leq C, \quad \forall\, q\in (1,+\infty).\label{intvw2}
 \eea
 Moreover, for any $T>0$,
 \be
   \|\rho(t)\|_{W^{1,q}}\leq e^{C T}\|\rho_0\|_{W^{1,q}}, \quad \forall\, t\in [0,T].\label{esrhof}
 \ee
  The constant $C$ in the above estimates depends on $\nu, \|v_0\|_{\mathbf{H}^2}, \|d_0
 \|_{\mathbf{H}^3}, \eta, \underline{\rho}, \bar{\rho}$ and also $\|\phi\|_{H^2}$ (under (F1)) or $\|\mathbf{g}\|_{L^2(0,+\infty;\, \mathbf{H}^1)}$, $\|\mathbf{g}_t\|_{L^2(0,+\infty;\, \mathbf{L}^2)}$ (under (F2)), but it is independent of $t$.
   \ep
  \begin{proof}
  It follows from \eqref{intbd1}, the definition of $B(t)$ and \eqref{eqv2} that
  \be
  \sup_{t\geq 0} \int_t^{t+1} B(\tau) d\tau \leq C,\label{intB1}
  \ee
  where $C$ is a constant depending on $\nu, \|v_0\|_{\mathbf{H}^1}, \|d_0
 \|_{\mathbf{H}^2}, \eta, \underline{\rho}, \bar{\rho}$ and also $\|\phi\|_{H^1}$ (under (F1)) or $\|\mathbf{g}\|_{L^2(0,+\infty;\, \mathbf{L}^2)}$ (under (F2)).
  Applying the uniform Gronwall lemma again, we infer from \eqref{eB1} or \eqref{eB2}, and \eqref{intB1} that
 \be
 B(t)\leq C, \quad \forall\,t\geq 0,
 \ee
 where $C$ is a constant depending on $\nu, \|v_0\|_{\mathbf{H}^2}, \|d_0
 \|_{\mathbf{H}^3}, \eta, \underline{\rho}, \bar{\rho}$, and also $\|\phi\|_{H^2}$ (under (F1)) or $\|\mathbf{g}\|_{L^2(0,+\infty;\, \mathbf{H}^1)}$, $\|\mathbf{g}_t\|_{L^2(0,+\infty;\, \mathbf{L}^2)}$ (under (F2)). This yields the estimate \eqref{ubdd1}. Integrating \eqref{eB1} or \eqref{eB2} from $t$ to $t+1$, we conclude \eqref{intbd2}.
 Recalling \eqref{vvv}, we deduce from Lemma \ref{S} and the estimates \eqref{ubdd},  \eqref{eqv1} that
 \bea
 \|v\|_{\mathbf{H}^2}&\leq& C\left(\|\rho v_t\|+\|\rho(v\cdot\nabla v)\|+\left\|(\Delta d-f(d))\cdot \nabla d\right\|+\|\rho\mathbf{g}\|\right)\non\\
 &\leq& C\bar{\rho}(\|v_t\|+ \|v\|_{\mathbf{L}^4}\|\nabla v\|_{\mathbf{L}^4}+\|\mathbf{g}\|)+C \|\nabla d\|_{\mathbf{L}^4}\|\Delta d-f(d)\|_{\mathbf{L}^4}\non\\
 &\leq&
 C(\|v_t\|+\|\nabla v\|\|v\|_{\mathbf{H}^2}^\frac12\|\nabla v\|^\frac12+\|\mathbf{g}\|)\non\\
 && +C\|d\|_{\mathbf{H}^2}(\|\nabla(\Delta d-f(d))\|^\frac12\|\Delta d-f(d)\|^\frac12) \non\\
 &\leq & \frac12\|v\|_{\mathbf{H}^2}+C(\|v_t\|+ \|\nabla d_t\|+\|\nabla v\|+\|\mathbf{g}\|).\label{vh2b}
 \eea
Using \eqref{ubdd1} and the facts that either under (F1) $\|\mathbf{g}\|\leq \|\phi\|_{H^1}$, or under (F2) $\|\mathbf{g}\|$ is bounded by a constant depending on $\|\mathbf{g}\|_{L^2(0,+\infty;\, \mathbf{L}^2)}$ and $\|\mathbf{g}_t\|_{L^2(0,+\infty;\, \mathbf{L}^2)}$, we get the uniform estimate for $\|v\|_{\mathbf{H}^2}$.

 Next, by the elliptic estimate, we infer from \eqref{dtd} that
 \bea
 \|\nabla d\|_{\mathbf{H}^2}&\leq& C\left(\|\nabla d_t\|+\|\nabla (v\cdot \nabla d)\|+\|f'(d)\nabla d\|+ \|\nabla d_0\|_{\mathbf{H}^\frac32(\Gamma)}+\|\nabla d\|\right)\non\\
 &\leq& C+ C\|\nabla v\|_{\mathbf{L}^4}\|\nabla d\|_{\mathbf{L}^4}+C\|v\|_{\mathbf{L}^\infty}\|d\|_{\mathbf{H}^2}\non\\
 &\leq& C+ C\|v\|_{\mathbf{H}^2}\|d\|_{\mathbf{H}^2}.\label{dh3b}
 \eea
 Combining \eqref{vh2b} and \eqref{dh3b}, we obtain the uniform estimate \eqref{ubdd2}. In a similar manner to \eqref{vh2b}, for $q\in (1, +\infty)$, using         \eqref{ubdd2}, we get
 \bea
 \|v\|_{\mathbf{W}^{2,q}}&\leq&  C\left(\|\rho v_t\|_{\mathbf{L}^q} +\|\rho(v\cdot\nabla v)\|_{\mathbf{L}^q}+\left\|(\Delta d-f(d))\cdot\nabla d\right\|_{\mathbf{L}^q}+\|\rho\mathbf{g}\|_{\mathbf{L}^q}\right)\non\\
 &\leq& C\bar{\rho}(\|v_t\|_{\mathbf{L}^q}+ \|v\|_{\mathbf{L}^\frac{q}{2}}\|\nabla v\|_{\mathbf{L}^\frac{q}{2}}+\|\mathbf{g}\|_{\mathbf{L}^q})+C \|\nabla d\|_{\mathbf{L}^\frac{q}{2}}\|\Delta d-f(d)\|_{\mathbf{L}^\frac{q}{2}}\non\\
 &\leq&
 C(\|\nabla v_t\|+\|v\|_{\mathbf{H}^2}+\|\nabla (\Delta d-f(d))\|+\|\mathbf{g}\|_{\mathbf{L}^q}) \non\\
 &\leq&C (\|\nabla v_t\|+\|\nabla (\Delta d-f(d))\|+\|\mathbf{g}\|_{\mathbf{H}^1}+\|\nabla v\| ).
 \label{vw2q}
 \eea
 Thus, \eqref{intvw2} is a consequence of \eqref{vw2q}, \eqref{intbd1} and \eqref{intbd2}, where the bound either depends on $\|\phi\|_{H^2}$ (under (F1)) or $\|\mathbf{g}\|_{L^2(0,+\infty;\, \mathbf{H}^1)}$ (under (F2)). Finally, for any $T>0$,
 \bea
 \int_0^t \|\nabla v(\tau)\|_{\mathbf{L}^\infty} d\tau &\leq& C\int_0^T \|v(\tau)\|_{\mathbf{W}^{2,r}}d\tau\leq CT^\frac12 \left(\int_0^T \|v(\tau)\|_{\mathbf{W}^{2,r}}^2d\tau\right)^\frac12\non\\
 &\leq& CT, \quad \forall \, t\in [0,T], \ \text{for some}\ r>2.\non
 \eea
 Then we infer from \eqref{w1qrho} that \eqref{esrhof} holds.  The proof is complete.
  \end{proof}
  \subsection{Proof of Theorem \ref{main2d}} Under the assumptions in Theorem \ref{main2d}, we can derive estimates for the approximate solutions to problem \eqref{0}--\eqref{5} $(\rho_m, v_m, d_m)$ as in Propositions \ref{v0h1}, \ref{es12} and \ref{est11}, which are independent of the parameter $m$. By extracting a subsequence and passing to limit as $m\to +\infty$, we can obtain a global strong solution $(\rho, v, d)$ to problem \eqref{0}--\eqref{5} in a standard way such that
   \bea
  && \rho \in L^\infty([0,T], W^{1,r}(\Omega)), \non\\
  && v\in L^\infty([0,T]; \mathbf{H}^2(\Omega)\cap V)\cap L^2(0,T; \mathbf{W}^{2,q}), \quad v_t\in L^2(0,T; \mathbf{V}), \ q\in (1, +\infty),\non\\
  && d\in L^\infty([0,T]; \mathbf{H}^3(\Omega)), \quad d_t\in L^\infty([0,T]; \mathbf{H}^1_0(\Omega))\cap L^2(0,T; \mathbf{H}^2),\non\\
  && 0<\underline{\rho}\leq \rho(x,t )\leq \bar{\rho},\quad |d(x,t)|\leq 1,\quad \forall\,
 (x,t)\in \Omega\times[0,T].\non
 \eea
 Due to Lemma \ref{trans}, we infer that $\rho \in C([0,T], W^{1,r}(\Omega))$. It is not difficult to see from \eqref{ubdd2} and \eqref{intbd2} that
 $d\in L^2(0,T; \mathbf{H}^4)$, which together with $d_t\in L^2(0,T; \mathbf{H}^2) $ yields that $d\in C([0,T]; \mathbf{H}^3)$. By the regularity of the Stokes operator,
  \bea
 \|v\|_{\mathbf{H}^{3}}&\leq&  C\left(\|\rho v_t\|_{\mathbf{H}^1} +\|\rho(v\cdot\nabla v)\|_{\mathbf{H}^1}+\left\|(\Delta d-f(d))\cdot \nabla d\right\|_{\mathbf{H}^1}+\|\rho\mathbf{g}\|_{\mathbf{H}^1}\right)\non\\
 &\leq& C\bar{\rho}(\|v_t\|_{\mathbf{H}^1}+ \|v\cdot\nabla v\|_{\mathbf{H}^1}+\|\mathbf{g}\|_{\mathbf{H}^1})\non\\
 && + C\|\nabla \rho\|_{\mathbf{L}^{q'}}(\|v_t\|_{\mathbf{L}^q}+ \|v\cdot\nabla v\|_{\mathbf{L}^q}+\|\mathbf{g}\|_{\mathbf{L}^q})+C \|\nabla d\|_{\mathbf{L}^\infty}\|\Delta d-f(d)\|_{\mathbf{H}^1}\non\\
 && +C\|\nabla^2 d\|_{\mathbf{L}^{q'}}\|\Delta d-f(d)\|_{\mathbf{L}^q}, \ \ \text{for some}\ 2<q,q'<+\infty, \ \frac{1}{q}+\frac{1}{q'}=\frac12 \non
 \eea
and estimates \eqref{ubdd2}, \eqref{intbd2}, \eqref{esrhof}, we infer that $v\in L^2(0,T; \mathbf{H}^3)$. This fact and $v_t\in L^2(0,T; \mathbf{V})$ yields the continuity $v\in C([0,T], \mathbf{H}^2)$.

Finally, we briefly show that the strong solution is indeed unique. Suppose $(\rho, v, d)$ and $(\tilde{\rho}, \tilde{v}, \tilde{d})$ are two strong solutions corresponding to the same initial data $(\rho_0, v_0, d_0)$. One can easily check that all the computations in \cite[Section 4]{DQS} can be verified due to the regularity of the two strong solutions. Denote $\delta \rho=\rho-\tilde{\rho}$, $\delta v=v-\tilde{v}$ and $\delta d=d-\tilde{d}$. We have for any $t\in [0,T]$ (cf. e.g., \cite[(4.46)]{DQS}):
 \bea
 &&\frac12\|\delta \rho(t)\|^2+\frac12\int_\Omega \tilde{\rho}(t)|\delta v(t)|^2dx +\frac12\|\nabla (\delta  d(t))\|^2\non\\
 &\leq& -\int_0^t\|\nabla (\delta v)\|^2dt-\int_0^t\|\Delta (\delta d)\|^2dt+\int_0^t\int_\Omega (\delta \rho) (\delta v)\nabla \rho dxdt\non\\
 && -\int_0^t\int_\Omega \tilde{\rho}\nabla v|\delta v|^2 dxdt  +\int_0^t\int_\Omega (\delta \rho)(\delta v)(v_t+v\cdot \nabla v)dxdt\non\\
 && +\int_0^t\int_\Omega v(\delta d)\Delta(\delta d) dxdt-\int_0^t\int_\Omega (\delta v)\nabla (\delta d) \Delta d dxdt\non\\
 && +\int_0^t\int_\Omega \Delta (\delta d)(f(d)-f(\tilde{d}))dxdt. \non
 \eea
 Using the H\"older inequality, Young inequality and the estimates \eqref{ubdd1}, \eqref{ubdd2}, \eqref{esrhof}, we easily get
 \be
 \|\delta \rho(t)\|^2+\int_\Omega \tilde{\rho}(t)|\delta v(t)|^2dx +\|\nabla (\delta  d(t))\|^2\leq C_T\int_0^t(\|\delta \rho\|^2+\|\delta v\|^2 +\|\nabla (\delta  d)\|^2)dt.\non
 \ee
 Due to the positivity of the density $\tilde{\rho}\geq \underline{\rho}>0$, the uniqueness follows from the above estimate and the Gronwall inequality. The proof is complete. \quad $\square$

\section{Long-time behavior}
\setcounter{equation}{0}
In this section, we study the long-time behavior of global strong solutions to problem \eqref{0}--\eqref{5}.

\bp \label{comconv1}
Under the assumptions in Theorem \ref{main2d}, either $\mathbf{g}$ satisfies (F1) or (F2), the global strong solution to \eqref{0}--\eqref{5} has the following decay property:
  \be
  \lim_{t\rightarrow +\infty} (\|v(t)\|_{\mathbf{H}^1}+ \|\Delta d(t)-f(d(t))\|)=0. \label{con1}
  \ee
\ep
\begin{proof}
For both cases (F1) and (F2), \eqref{intA} implies that $\int_0^{+\infty} A(t) dt<+\infty$. Then it follows from \eqref{he} and \cite[Lemma 6.2.1]{Z04} that $\lim_{t\to+\infty}A(t)=0$, which together with the Poincar\'e inequality yields \eqref{con1}.
\end{proof}

\bp \label{comconv2}
Under the assumptions in Theorem \ref{main2d}, either $\mathbf{g}$ satisfies (F1) or (F2), the global strong solution to \eqref{0}--\eqref{5} has the following decay property:
  \be
  \lim_{t\rightarrow +\infty} (\|v_t(t)\|_{\mathbf{H}^1}+ \|\nabla(\Delta d(t)-f(d(t)))\|+\|d_t\|_{\mathbf{H}^1})=0. \label{con2}
  \ee
\ep
 \begin{proof}
 First, we look at the case that $\mathbf{g}$ satisfies (F2). We have already proved that $A(t)$ is bounded for all time and $A(t)\in L^1(0, +\infty)$. Therefore, integrating \eqref{he} from $0$ to $+\infty$ with respect to time, we infer from (F2) that
 \bea  && \int_0^{+\infty} (\|\rho^\frac12 v_t(t)\|^2+\gamma\|\nabla(\Delta d(t)-f(d(t)))\|^2)dt \non\\
 &\leq& A(0)+C\left(\sup_{t\geq 0}A(t)+1\right)\int_0^{+\infty} A(t) dt+C\int_0^{+\infty} \|\mathbf{g}(t) \|^2dt<+\infty.
 \eea
 Recalling \eqref{eqv2}, we obtain
 \be \int_0^{+\infty}B(t) dt<+\infty. \label{intB}
 \ee
 Then using \eqref{eB2}, assumption (F2) and \cite[Lemma 6.2.1]{Z04}, we conclude that
 \be
 \lim_{t\to +\infty} B(t)=0,
 \ee
 which combined with Poincar\'e inequality and \eqref{eqv1}, \eqref{eqv2} yields \eqref{con2}.

 Next, we deal with the case when $\mathbf{g}$ satisfies (F1). We only need to show that \eqref{intB} still holds in this case.
 For this purpose, we re-estimate the term $I_2$, $I_3$, $I_4$ and $I_5$ in \eqref{dA} using the higher-order estimate \eqref{ubdd2} instead of the lower-order one \eqref{ubdd}.
 As in \cite{ZK}, we infer from the transport equation \eqref{0} and integration by parts that
 \bea
 I_3&=& \frac{d}{dt}\int_\Omega \rho\nabla \phi \cdot v dx-\int_\Omega \rho_t\nabla \phi \cdot v dx\non\\
&=& \frac{d}{dt}\int_\Omega \rho\nabla \phi \cdot v dx - \int_\Omega \rho v\cdot \nabla (\nabla \phi \cdot v)dx\non\\
 &\leq& \frac{d}{dt}\int_\Omega \rho\nabla \phi \cdot v dx+ \bar{\rho}(\|v\|_{\mathbf{L}^4}^2\|\phi\|_{H^2}+\|v\|_{\mathbf{L}^4}\|\nabla \phi\|_{\mathbf{L}^4}\|\nabla v\|)\non\\
 &\leq& \frac{d}{dt}\int_\Omega \rho\nabla \phi \cdot v dx+ C(\|v\|\|\nabla v\|+\|v\|^\frac12\|\nabla v\|^\frac32)\non\\
 &\leq& \frac{d}{dt}\int_\Omega \rho\nabla \phi \cdot v dx+C\|\nabla v\|^2.\non
 \eea
 Next, since $\Delta d-f(d)|_\Gamma=0$, we obtain that
 \bea
 I_2+I_4+I_5
 &=& - \lambda \int_\Omega [ (\Delta d-f(d))\cdot \nabla d]\cdot  v_tdx -\int_\Omega \Delta ( v\cdot \nabla d)\cdot (\Delta d-f(d)) dx\non\\
 &=& - \lambda \int_\Omega [(\Delta d-f(d))\cdot \nabla d]\cdot  v_tdx +\int_\Omega \nabla  ( v\cdot \nabla d)\cdot \nabla (\Delta d-f(d)) dx \non\\
 &=&   - \lambda \int_\Omega [(\Delta d-f(d))\cdot \nabla d]\cdot  v_tdx +\int_\Omega \nabla_k v_i \nabla_i d_j \nabla_k (\Delta d_j-f_j(d)) dx\non\\
 && +  \int_\Omega    v_i\nabla_k \nabla_i d_j \nabla_k (\Delta d_j-f_j(d)) dx\non\\
 &:=& I'_2+I'_4+I'_5,
 \eea
 where
 \bea
 I'_2&\leq&   C\underline{\rho}^\frac12 \|\rho^\frac12 v_t\|\|(\Delta d-f(d))\cdot \nabla d\|\non\\
 &\leq&  C\|\rho^\frac12 v_t\|(\|\Delta d-f(d)\|^\frac12+1)\|\Delta d-f(d)\|^\frac12\|\nabla (\Delta d-f(d))\|^\frac12\non\\
  &\leq& \frac38\int_\Omega \rho|v_t|^2dx+ \frac{\gamma}{8} \|\nabla (\Delta d-f(d))\|^2+C\|\Delta d-f(d)\|^4+C\|\Delta d-f(d)\|^2,\non
 \eea
 \bea
 I'_4+I'_5&\leq& \|\nabla v\|\|\nabla d\|_{\mathbf{L}^\infty}\|\nabla (\Delta d-f(d))\|+\| v\|_{\mathbf{L}^4}\|\nabla ^2 d\|_{\mathbf{L}^4}\|\nabla (\Delta d-f(d))\|\non\\
 &\leq& C\|\nabla v\|\|d\|_{\mathbf{H}^3}\|\nabla (\Delta d-f(d))\|\non\\
 &\leq& \frac{\gamma}{4} \|\nabla (\Delta d-f(d))\|^2+C\|\nabla v\|^2,\non
 \eea
 Replacing the original estimates for $I_2$,..., $I_5$ by $I_2'$,..., $I_5'$, we arrive at the following inequality
  \bea && \frac{d}{dt}\left(A(t)-2\int_\Omega \rho\nabla \phi\cdot v dx\right)+\|\rho^\frac12 v_t(t)\|^2+\gamma\|\nabla(\Delta d(t)-f(d(t)))\|^2\non\\
  &\leq& C(A^2(t)+A(t)), \label{hea}
 \eea
  where $C$ is a constant depending on $\|v_0\|, \|d_0 \|_{\mathbf{H}^1}, \eta, \underline{\rho}, \bar{\rho}, \nu, \Omega$ and also $\|\phi\|_{H^2}$. Integrating \eqref{hea} with respect to time from $0$ to $+\infty$, we deduce that
 \bea
&& \int_0^{+\infty}\|\rho^\frac12 v_t(t)\|^2+\gamma\|\nabla(\Delta d(t)-f(d(t)))\|^2dt \non\\
 &\leq& A(0)-2\int_\Omega \rho_0\nabla \phi \cdot v_0 dx+ 2\left|\int_\Omega \rho\nabla \phi v dx\right|\non\\
 && +C\left(\sup_{t\geq 0}A(t)+1\right)\int_0^{+\infty} A(t) dt\non\\
 &<&+\infty,\non
 \eea
 which again implies \eqref{intB}. Using the same argument as for the previous case, we conclude the decay property  \eqref{con2}. The proof is complete.
 \end{proof}

 \subsection{Proof of Theorem \ref{conv1}}
 According to Propositions \ref{comconv1} and \ref{comconv2}, it remains to show the convergence for the density function $\rho$  and the director vector $d$.

 As a direct consequence of the uniform-in-time estimate \eqref{bdr} and weak compactness of bounded sets in $L^q$ ($1<q<+\infty$), we know that for any sequence $\{t_i\}\nearrow+\infty$, there is a subsequence $\{t'_i\}\nearrow+\infty$ such that $\rho(t'_i)$ weakly converge to a certain $\rho_\infty$ in $L^q$. On the other hand, due to Lemma \ref{trans}, $\|\rho(t'_i)\|_{L^q}=\|\rho_\infty\|_{L^q}=\|\rho_0\|_{L^q}$, since $L^q$ ($1<q<+\infty$) is a uniformly convex Banach space, we can conclude that  $\rho(t'_i)$ actually strongly converge to $\rho_\infty$ in $L^q$ (cf. \cite[Lemma 3.1.6]{Z04}), namely, \eqref{conrho1} holds true.

 The uniform-in-time estimate \eqref{ubdd2} yields that for any sequence $\{t_i\}\nearrow+\infty$, there is a subsequence $\{t'_i\}\nearrow+\infty$ such that $d(t'_i)$ strongly converge to a certain $d_\infty$ in $\mathbf{H}^2$. Then by \eqref{con1}, we see that
  \bea
  0&\leq& \|-\Delta d_\infty+f(d_\infty)\|\non\\
  &\leq&  \|-\Delta d(t'_i)+f(d(t'_i))\|+\|\Delta d(t'_i)-\Delta d_\infty\|+\|f(d(t'_i))-f(d_\infty)\|\non\\
  &\leq&  \|-\Delta d(t'_i)+f(d(t'_i))\|+C\|d(t'_i)-d_\infty\|_{\mathbf{H}^2}\non\\
  &\to& 0, \quad \text{as}\ \ t'_i\nearrow+\infty.\non
  \eea
  Thus, $d_\infty$ satisfies the stationary problem \eqref{staa}. The $\mathbf{H}^3$ convergence follows from \eqref{con2} and the fact that
  \bea
  \|\nabla \Delta (d(t'_i)-d_\infty)\|&\leq&  \|\nabla (-\Delta d(t'_i)+f(d(t'_i)))\|+ \|\nabla (f(d(t'_i))-f(d_\infty))\|\non\\
  &\leq& \|\nabla (-\Delta d(t'_i)+f(d(t'_i)))\|+ C\|d(t'_i)-d_\infty\|_{\mathbf{H}^2}\non\\
   &\to& 0, \quad \text{as}\ \ t'_i\nearrow+\infty.\non
  \eea
  The proof of Theorem \ref{conv1} is complete.
  \begin{remark} \label{the2re}
  (1) When the external force $\mathbf{g}$ is a gradient field independent of time, our results show that the velocity field and its time derivative will converge to zero as time goes to infinity, this coincides with the result for the $2D$ density-dependent incompressible Navier--Stokes equations \cite{ZK}.

  (2) For the density function $\rho$, we are only able to obtain some partial results that it will sequentially converge to a certain function in $L^q$ norm (differential sequences may have different limit points). One possible sufficient condition for the convergence in time is that $\|\rho_t\|_{H^{-1}}\in L^1(0,+\infty)$. On the other hand, we observe that for any function $\varphi\in H^1_0(\Omega)$,
 $$ \int_\Omega \rho_t \varphi dx= -\int_\Omega \nabla \cdot(\rho v) \varphi dx= \int_\Omega \rho v\cdot \nabla \varphi dx \leq \bar{\rho}\|v\|\|\nabla \varphi\|,$$
 which yields $\|\rho_t\|_{H^{-1}}\leq C\|v\|$.
 Hence, if we can prove $\|v\|\in L^1(0,+\infty)$, then we get the convergence of $\rho$ as time tends to infinity in the space $H^{-1}$ and thus in $L^q$ due to the uniqueness of limit. However, we do not know this $L^1$-integrability condition on $v$ from the above proof.

 (3) For the director vector $d$, we are able to show the decay of its time derivative and sequential convergence of itself to a steady state that is a solution to the stationary problem \eqref{staa}. If we know that problem \eqref{staa} admits a unique solution, then $d$ will converge to it as time goes to infinity. However, in general we cannot expect the uniqueness of solutions to the stationary Ginzburg--Landau equation \eqref{staa}, unless, for instance, the parameter $\eta$ in $f(d)$ is sufficiently large.
  \end{remark}

 \subsection{Proof of Theorem \ref{conv}}
 In this subsection, we provide the proof for Theorem \ref{conv}. Now the external force $\mathbf{g}$ is an asymptotically autonomous one such that it satisfies (F2) and \eqref{ggg}. Different from the previous case with time-independent force, we shall see that one can obtain convergence for the density function and director vector. The proof is based on an appropriate generalization of the {\L}ojasiewicz-Simon approach for gradient-like systems (cf. e.g., \cite{HT01}).

Denote
 \be
 E(d)=\frac12\|\nabla d\|^2 + \int_\Omega F(d)dx.\label{EDD}
 \ee
 It is straightforward to check that any solution to the stationary problem \eqref{staa} is a
 critical point of the energy functional $E(d)$, and conversely, the
 critical point of $E(d)$ is a solution to \eqref{staa} (cf. \cite{W10}).
 Besides, regularity of solutions to \eqref{staa} has been shown in \cite{LL95}
 such that $d$ is smooth on $\Omega$ provided that $d_0$ is smooth on $\Gamma$. Below we shall make use the following \L ojasiewicz--Simon type inequality (cf. \cite{W10}):
\bl \label{ls}
 Let $\psi$ be a critical point of $E(d)$. There exist constants
 $\theta\in(0,\frac12)$ and $\beta>0$ depending on $\psi$ such that
 for any $d\in \mathbf{H}^2(\Omega)$ satisfying $d|_\Gamma=d_0(x)$ and $\|d-\psi\|_{\mathbf{H}^2}<\beta$,
 there holds
 \be
 \|-\Delta d+f(d)\|\geq
 |E(d)-E(\psi)|^{1-\theta}.\label{LSE}
 \ee
 \el
The uniform estimate \eqref{ubdd2} yields that the set of all limit points of the trajectory $\{ d(t): t\geq 0\}$ denoted by $\omega(d):=\bigcap_{s\geq 0}\overline{\{d(t)\in \mathbf{H}^2(\Omega): t\geq s\}}$ is non-empty and compact in $\mathbf{H}^2(\Omega)$. Besides, similar to Section 4.1, we can see that each element in $\omega(d)$ is a solution to problem \eqref{staa} and thus is a critical point of $E(d)$. We infer from \eqref{lya1} that
\be
\mathcal{E}(t_1)-\mathcal{E}(t_2)\leq C\int_{t_2}^{t_1} \|\mathbf{g}\|^2 dt, \quad \forall\, t_1>t_2>0.
\ee
Since $\|\mathbf{g}\|\in L^2(0,+\infty)$, it follows that $\mathcal{E}(t)$ converges to a certain constant $\mathcal{E}_\infty$ as $t\to+\infty$. Recalling the convergence of velocity field \eqref{con1}, we deduce that $E(d)$ is indeed a constant on $\omega(d)$ such that $E(\psi)=\lambda^{-1}\mathcal{E}_\infty$, for all $\psi\in \omega(d)$.
Due to Lemma \ref{ls}, for every $\psi\in \omega(d)$, there exist some $\beta_\psi$ and $\theta_\psi\in (0, \frac12)$ that may depend on $\psi$ such that the inequality \eqref{LSE} holds for $d\in \mathbf{B}_{\beta_\psi}(\psi):=\{d\in \mathbf{H}^2(\Omega): d|_\Gamma=d_0, \|d-\psi\|_{\mathbf{H}^2}<\beta_\psi\}$ and $|E(d)-E(\psi)|\leq 1$. The union of balls $\{ \mathcal{B}_{\beta_\psi}(\psi): \psi\in \omega(d)\}$ forms an open cover of $\omega(d)$ and due to the compactness of $\omega(d)$, we can find a finite sub-cover $\{\mathbf{B}_{\beta_i}(\psi_i): i=1,2,...,m\}$ where the constants $\beta, \theta$ corresponding to $\psi_i$ in Lemma \ref{ls} are indexed by $i$. From the definition of $\omega(d)$, we know that there exist a sufficient large $t_0$ such that $d(t)\in \mathcal{U}:=\cup_{i=1}^m\mathcal{B}_{\beta_i}(\psi_i)$ for $t\geq t_0$. Taking $ \theta=\min_{i=1}^m\{\theta_i\}\in (0, \frac12)$, we get for all $t\geq t_0$
\be
 \|-\Delta d(t)+f(d(t))\|\geq
 |E(d(t))-\lambda^{-1}\mathcal{E}_\infty|^{1-\theta}.\label{LSEa}
 \ee
Denote $$\mathcal{Y}(t)^2=\frac{\nu}{2}\|\nabla v(t)\|^2+\lambda\gamma\|\Delta
 d(t)-f(d(t))\|^2,\quad z(t)=\int_t^{+\infty}\|\mathbf{g}(\tau)\|^2d\tau.$$
 Assumption \eqref{ggg} implies that $z(t)\leq C(1+t)^{-(1+\xi)}$ for $t\geq 0$.
Then by the basic energy inequality  \eqref{lya1}, we get
 \be
\mathcal{E}(t)-\mathcal{E}_\infty\geq \int_t^{+\infty} \mathcal{Y}(\tau)^2 d\tau-\frac{C_P^2\bar{\rho}^2}{2\nu}z(t)\geq  \int_t^{+\infty} \mathcal{Y}(\tau)^2 d\tau-C(1+t)^{-(1+\xi)}.\label{LSEc}
 \ee
On the other hand, using \eqref{LSEa}, the uniform estimates \eqref{ubdd} and the fact $\frac{1}{1-\theta}<2$, we obtain that
\bea
|\mathcal{E}(t)-\mathcal{E}_\infty|&\leq &
\frac12\|v\|^2+\lambda|E(d)-\lambda^{-1}\mathcal{E}_\infty|\non\\
&\leq& C\|\nabla v\|^{\frac{1}{1-\theta}}+ \lambda\|-\Delta d(t)+f(d(t))\|^{\frac{1}{1-\theta}}\non\\
&\leq& C\mathcal{Y}(t)^{\frac{1}{1-\theta}}, \quad \forall\, t\geq t_0.\label{LSEb}
\eea
Take $$\zeta=\min\left\{\theta, \frac{\xi}{2(1+\xi)}\right\}\in \left(0, \frac12\right).$$
It is easy to check that
\be
\int_t^{+\infty} (1+\tau)^{-2(1+\xi)(1-\zeta)}d\tau\leq \int_t^{+\infty}(1+\tau)^{-(2+\xi)}d\tau \leq (1+t)^{-(1+\xi)}, \quad \forall\, t\geq 0.\label{LSEd}
\ee
We now set
\be Z(t)=\mathcal{Y}(t)+(1+t)^{-(1+\xi)(1-\zeta)}.\non\ee
Since  $\lim_{t\to+\infty} \mathcal{Y}(t)=0$ (see  \eqref{con1}), it follows from \eqref{LSEc}, \eqref{LSEb} and \eqref{LSEd} that
\bea
\int_t^{+\infty}Z(\tau)^2d\tau &\leq& C\mathcal{Y}(t)^\frac{1}{1-\theta}+C(1+t)^{-(1+\xi)}\leq C\mathcal{Y}(t)^\frac{1}{1-\zeta}+C(1+t)^{-(1+\xi)}\non\\
&\leq& CZ(t)^\frac{1}{1-\zeta}, \quad \forall\, t\geq t_0.\label{Z}
\eea
Recall the following result (cf. \cite[Lemma 4.1]{HT01} or \cite[Lemma 7.1]{FS})
\bl\label{AEle}
Let $\zeta\in (0,\frac12)$. Assume that $Z\geq 0$ be a measurable function on $(0,+\infty)$, $Z\in L^2(\mathbb{R}^+)$ and there exist $C>0$ and $t_0\geq 0$ such that
 \be
\int_t^{+\infty}Z(\tau)^2 d\tau\leq C Z(t)^\frac{1}{1-\zeta},\quad \text{for a.e.}\ \  t\geq t_0.\non\ee
 Then $Z\in L^1(t_0, +\infty)$.
\el
\noindent We conclude from \eqref{Z} and Lemma \ref{AEle} that
$$\int_{t_0}^{+\infty}Z(t)<+\infty.$$
 Since $\xi>0$, it holds
 \be
\int_{t_0}^{+\infty} (1+t)^{-(1+\xi)(1-\zeta)}dt\leq \int_{t_0}^{+\infty}(1+t)^{-\frac12(2+\xi)}dt \leq 2(1+t_0)^{-(1+\xi)}<+\infty,\ \text{for}\ t_0>0,\non
 \ee
which implies that
\be
\int_{t_0}^{+\infty} \|\nabla v(t)\|+\|\Delta d(t)-f(d(t))\| dt<+\infty.\label{L1}
\ee
On the other hand, it follows from equation \eqref{2} that
 \be \|d_t\|\leq C(\|v\|_{\mathbf{L}^4}\|\nabla d\|_{\mathbf{L}^4}+\|-\Delta d+f(d)\|)\leq  C(\|\nabla v\|+ \|-\Delta d+f(d)\|).\label{dt}
 \ee
As a consequence,
 \be \int_{t_0}^\infty \|d_t(t)\| dt<+\infty,
 \ee
which easily implies that as $t\rightarrow +\infty$, $d(t)$
converges strongly in $\mathbf{L}^2(\Omega)$.  By compactness of $d(t)$ in $\mathbf{H}^2(\Omega)$, we deduce that
 \be \lim_{t\rightarrow +\infty}
 \|d(t)-d_\infty\|_{\mathbf{H}^2}=0,\label{vcon1}
 \ee
 where $d_\infty$ is a solution to problem \eqref{staa}. Recalling the uniform  estimate \eqref{ubdd2}, it holds
 \bea \|\nabla \Delta d- \nabla \Delta d_\infty\|&\leq& \| \nabla (\Delta d- \Delta d_\infty
 -f(d)+f(d_\infty))\|+ \|\nabla (f(d)-f(d_\infty))\|\non
 \\
 &\leq& \|\nabla (\Delta
 d-f(d))\|+C\|d-d_\infty\|_{\mathbf{H}^2}.\non
 \eea
 The above estimate together with \eqref{con2} and \eqref{vcon1} yields
 \be \lim_{t\rightarrow +\infty}
 \|d(t)-d_\infty\|_{\mathbf{H}^3}=0.\label{cond3}
 \ee
 It follows from Remark \ref{the2re} and  \eqref{L1} that $\|\rho_t\|_{H^{-1}}\in L^1(t_0,+\infty)$. Thus $\rho(t)$ converges strongly in $H^{-1}$ as $t\to +\infty$. By an argument similar to that in the proof of Theorem \ref{conv1}, we conclude that
  \be \lim_{t\rightarrow +\infty}
 \|\rho(t)-\rho_\infty\|_{L^q}=0, \quad q\in (1,+\infty).
  \ee

 Next, we prove the convergence rate. Denote
 \be
 \mathcal{K}(t)= \mathcal{E}(t)-\mathcal{E}_\infty+\frac{C_P^2\bar{\rho}^2}{2\nu}\int_t^{+\infty}\|\mathbf{g}(\tau)\|^2d\tau. \non
 \ee
 It follows from the basic energy inequality \eqref{lya1} that
 \be
 \frac{d}{dt}\mathcal{K}(t)+\mathcal{Y}(t)^2\leq 0.\label{YY}
 \ee
 Thus, $\mathcal{K}(t)$ is decreasing on $[0,+\infty)$ and $\mathcal{K}(t)\to 0$ as $t\to +\infty$. Recalling the definition of $t_0$, for $t\geq t_0$, we deduce from \eqref{ggg} and \eqref{LSEb} that
 \bea
 \mathcal{K}(t)^{2(1-\zeta)}&\leq& C\mathcal{Y}(t)^\frac{2(1-\zeta)}{1-\theta}+C(1+t)^{-2(1-\zeta)(1+\xi)}\non\\
 &\leq& C\mathcal{Y}(t)^2+C(1+t)^{-2(1-\zeta)(1+\xi)}\non\\
 &\leq& -C\frac{d}{dt}\mathcal{K}(t)+C(1+t)^{-2(1-\zeta)(1+\xi)},\label{KKK}
 \eea
 where we have used the fact that $\frac{2(1-\zeta)}{1-\theta}\geq 2$. It follows from the ordinary differential inequality \eqref{KKK} and \cite[Lemma 2.6]{Ben} that
 \be
 \mathcal{K}(t)\leq C(1+t)^{-\iota}, \quad \forall\, t\geq t_0,\label{rate1}
 \ee with the exponent given by
 \be
  \iota=\min\left\{\frac{1}{1-2\zeta}, 1+\xi\right\}=\min\left\{\frac{1}{1-2\theta}, 1+\xi\right\}.\non
 \ee
 We infer from \eqref{YY} that for any $t\geq t_0$,
 \be
 \int_t^{2t} \mathcal{Y}(\tau) d\tau \leq t^{\frac12}\left(\int_t^{2t} \mathcal{Y}(\tau)^2 d\tau\right)^\frac12 \leq  Ct^\frac12 \mathcal{K}(t)^\frac12\leq  C(1+t)^{\frac{1-\iota}{2}},\non
 \ee
 where $$\kappa=\frac{\iota-1}{2}=\min\left\{\frac{\theta}{1-2\theta}, \frac{\xi}{2}\right\}>0.$$
 It holds
 \be
 \int_t^{+\infty} \mathcal{Y}(\tau) d\tau\leq \sum_{j=0}^{+\infty} \int_{2^j t}^{2^{j+1}t} \mathcal{Y}(\tau) d\tau
 \leq C\sum_{j=0}^{+\infty} (2^j t)^{-\kappa}\leq  C(1+t)^{-\kappa}, \quad \forall\, t\geq t_0.\label{decayY}
 \ee
  Then by \eqref{dt}, we get
 \be
 \int_t^{+\infty} \|d_t(\tau)\|d \tau \leq \int_t^{+\infty} \mathcal{Y}(\tau) d\tau\leq  C(1+t)^{-\kappa}, \quad \forall\, t\geq t_0.\non\\
 \ee
 which together with \eqref{ubdd} yields the convergence rate of $d$ in $\mathbf{L}^2$
 \be
 \|d(t)-d_\infty\|\leq C(1+t)^{-\kappa}, \quad \forall\, t\geq 0.\label{rated1}
 \ee
  Besides, using the Poincar\'e inequality, we infer from Remark \ref{the2re} and \eqref{decayY} that
\be
\int_t^{+\infty}\|\rho_t(\tau)\|_{H^{-1}}\leq C\int_{t}^{+\infty} \|v(\tau)\|d\tau \leq C\int_{t}^{+\infty} \|\mathcal{Y}(\tau)\|d\tau\leq  C(1+t)^{-\kappa}, \quad \forall\,t\geq 0,\non
\ee
which implies
\be
\|\rho(t)-\rho_\infty\|_{H^{-1}}\leq C(1+t)^{-\kappa}, \quad \forall\,t\geq 0.\non
\ee

 Taking advantage of the lower-order convergence rate of the director $d$ \eqref{rated1}, we are able to obtain decay estimate for $v$ as well as higher-order convergence rate on $d$. For this purpose, we make use the fact  $-\Delta d_\infty+f(d_\infty)=0$ and test \eqref{3} by
 $\lambda(-\Delta (d-d_\infty) + f(d)-f(d_\infty)+d-d_\infty)$ to get
\bea
 && \frac{d}{dt}\left(\frac{\lambda}{2}\|\nabla (d-d_\infty)\|^2+\frac{\lambda}{2}\|d-d_\infty\|^2+\lambda\int_\Omega F(d)-F(d_\infty)-f(d_\infty)\cdot(d-d_\infty)dx\right)\non\\
 && +\lambda\gamma\|-\Delta d+f(d)\|^2+\lambda\gamma\|\nabla (d-d_\infty)\|^2+\lambda \int_\Omega (v\cdot \nabla d)\cdot \Delta d dx\non\\
 &=& -\lambda\int_\Omega (v\cdot \nabla d) \cdot (d-d_\infty)dx- \lambda\gamma\int_\Omega (f(d)-f(d_\infty))\cdot (d-d_\infty)dx.\label{b2a}
 \eea
 Using the uniform estimate \eqref{ubdd2} and Sobolev embeddings, we can estimate the right-hand side of \eqref{b2a} as follows
 \be
 \text{R.H.S of \eqref{b2a}}\leq \frac{\nu}{4}\|\nabla v\|^2+C\|d-d_\infty\|^2.\label{rhs}
 \ee
 Adding \eqref{b2a} with \eqref{b0} and \eqref{b1}, we infer from \eqref{rhs} that
 \be
 \frac{d}{dt}\mathcal{Q}(t)+C_1(A(t)+\|\nabla (d-d_\infty)\|^2)\leq C_2(\|d-d_\infty\|^2+\|\mathbf{g}\|^2), \label{dQ}
 \ee
 where
 \bea
 \mathcal{Q}(t)&=&\frac{1}{2}\int_\Omega \rho(t) |v(t)|^2 dx+ \frac{\lambda}{2}\|\nabla (d(t)-d_\infty)\|^2+\frac{\lambda}{2}\|d(t)-d_\infty\|^2\non\\
 &&+\lambda\int_\Omega F(d(t))-F(d_\infty)-f(d_\infty)\cdot(d(t)-d_\infty)dx.\non
 \eea
 By the uniform estimate \eqref{ubdd2} and the Taylor's formula, it is easy to see that $$\left|\int_\Omega F(d)-F(d_\infty)-f(d_\infty)(d-d_\infty)dx\right|\leq C_3\|d-d_\infty\|^2,$$
 which implies
 \be
 \mathcal{Q}(t)+C_3\|d(t)-d_\infty\|^2\geq \frac{1}{2}\int_\Omega \rho(t) |v(t)|^2 dx+ \frac{\lambda}{2}\|\nabla (d(t)-d_\infty)\|^2.\label{dQ1}
 \ee
 Since $A(t)$ is uniformly bounded in time, \eqref{he} can be rewritten as
 \be \frac{d}{dt} A(t)\leq C_4 A(t)+C_5 \|\mathbf{g}\|^2,\label{daa}\ee
 We deduce from \eqref{dQ}, \eqref{dQ1} and \eqref{daa} that
 \be
 \frac{d}{dt}\mathcal{Q}_1(t)+C_6 \mathcal{Q}_1(t)\leq C_7(\|d-d_\infty\|^2+\|\mathbf{g}\|^2),
 \ee
 where $\mathcal{Q}_1(t)=\mathcal{Q}(t)+\frac{C_1}{2C_4}A(t)$. As a consequence,
 \bea
 \mathcal{Q}_1(t)&\leq& e^{-C_6t}\left(\mathcal{Q}_1(0)+ C_7\int_0^t e^{C_6\tau}(\|d(\tau)-d_\infty\|^2+\|\mathbf{g}(\tau)\|^2)d\tau\right)\non\\
 &\leq& \mathcal{Q}_1(0) e^{-C_6t}+Ce^{-\frac{C_6t}{2}} \int_0^{\frac{t}{2}} (\|d(\tau)-d_\infty\|^2+\|\mathbf{g}(\tau)\|^2) d\tau \non\\
 && +C \left(\sup_{s\in[\frac{t}{2},t]} \|d(s)-d_\infty\|^2 \right) e^{-C_6t} \int_{\frac{t}{2}}^t e^{C_6\tau}  d\tau  +C\int_{\frac{t}{2}}^{+\infty} \|\mathbf{g}(t)\|^2 d\tau\non\\
 &\leq& \mathcal{Q}_1(0) e^{-C_6t}+e^{-\frac{C_6t}{2}}\left(C\int_0^\frac{t}{2} (1+\tau)^{-2\kappa} d\tau+C\right) \non\\&& +C\left(1+\frac{t}{2}\right)^{-2\kappa}+C \left(1+\frac{t}{2}\right)^{-(1+\xi)}\non\\
 &\leq& C(1+t)^{-2\kappa},\quad \forall\, t\geq 0.\label{decayQ1}
 \eea
 Recalling the definitions of $A(t)$, $\mathcal{Q}_1(t)$, $\mathcal{Q}(t)$, we conclude  from \eqref{decayQ1}, \eqref{dQ1}, \eqref{rated1} and \eqref{bdr} that
 \be
 \|v(t)\|_{\mathbf{H}^1}+\|d(t)-d_\infty\|_{\mathbf{H}^1}+\|\Delta d(t)-f(d(t))\|\leq C(1+t)^{-\kappa},\quad \forall\, t\geq 0.\label{rrr}
 \ee
 By the elliptic estimate
 \bea
 \|d-d_\infty\|_{\mathbf{H}^2}&\leq& C \|\Delta (d-d_\infty)\|\leq C\|\Delta d-f(d)\|+C\|f(d)-f(d_\infty)\|\non\\
 &\leq& C\|\Delta d-f(d)\|+\|d-d_\infty\|_{\mathbf{H}^1},
 \eea
 we get
 \be
 \|d(t)-d_\infty\|_{\mathbf{H}^2}\leq C(1+t)^{-\kappa},\quad \forall\, t\geq 0.
 \ee
 The proof of Theorem \ref{conv} is complete.

\begin{remark}
If the velocity field $v$ decays fast enough,
 we can obtain uniform-in-time $W^{1,r}$-estimate for $\rho$ $(1<r<+\infty)$.  For any $q\in (1, +\infty)$, we infer from \eqref{vw2q}, \eqref{intB} and (F2) that  $\|v(t)\|_{\mathbf{W}^{2,q}}\in L^2(0, +\infty)$.
If $\kappa>\frac32$, or in other words, $\xi>3$ and $\theta\in (\frac{3}{8}, \frac12)$, we just take $q\in (\frac{4\kappa-4}{2\kappa-3},+\infty)$ such that $\frac{(2q-4)\kappa}{3q-4}>1$. Then it follows from \eqref{rrr} that
 \bea
 && \int_0^{+\infty} \|\nabla v(t)\|_{\mathbf{L}^\infty} dt \non\\
 &\leq& C\int_0^{+\infty} \|v(t)\|_{\mathbf{W}^{2,q}}^\frac{q}{2(q-1)}\|\nabla v(t)\|^\frac{q-2}{2(q-1)} dt\non\\
 &\leq& C\left(\int_0^{+\infty} \|\nabla v(t)\|^\frac{2q-4}{3q-4}dt \right)^\frac{3q-4}{4(q-1)} \left(\int_0^{+\infty} \|v(t)\|_{\mathbf{W}^{2,q}}^2 dt \right)^\frac{q}{4(q-1)}\non\\
 &\leq& C\left(\int_0^{+\infty} (1+t)^{-\frac{(2q-4)\kappa}{3q-4}} dt \right)^\frac{3q-4}{4(q-1)}\non\\
 &\leq& C.\non
 \eea
Therefore, by \eqref{w1qrho}, we obtain that $\|\rho(t)\|_{W^{1,r}}\leq C\|\rho_0\|_{W^{1,r}}$ for $t\geq 0$.

\end{remark}

\section*{Acknowledgements} The authors would like to thank the referee for calling their attention to related works in the literature on this topic. X. Hu's research was partially supported by the National Science Foundation. H. Wu was
partially supported by National Science Foundation of China 11001058, SRFDP and ``Chen Guang" project supported by Shanghai
Municipal Education Commission and Shanghai Education Development
Foundation.


\end{document}